\documentclass[12pt]{amsart}
\usepackage{amsmath, amssymb, amscd}
\usepackage{graphicx}

\usepackage[all,import,rotate]{xy}
\usepackage{ragged2e}

\usepackage{mathpazo}
\usepackage{mathrsfs}
\renewcommand{\mathcal}{\mathscr}
\renewcommand{\coprod}{\bigcup}

\def\RCS$#1: #2 ${\expandafter\def\csname RCS#1\endcsname{#2}}
\RCS$Revision: 1.106 $
\RCS$Date: 2008/12/16 19:50:35 $
\renewcommand{\today}{\number\year /\number\month /\number\day}
\newcommand{\versioninfo}{Version \RCSRevision; Last commit \RCSDate; Compile \today }

\newcommand{\C}{{\mathbb C}}
\newcommand{\Q}{{\mathbb Q}}
\newcommand{\F}{{\mathcal F}}

\newcommand{\R}{{\mathbb R}}
\newcommand{\Z}{{\mathbb Z}}
\newcommand{\T}{{\mathcal T}}

\newcommand{\maps}{\colon\thinspace}
\DeclareMathOperator{\Hom}{Hom}
\DeclareMathOperator{\Comm}{Comm}
\DeclareMathOperator{\tr}{tr}

\DeclareMathOperator{\Isom}{Isom}

\DeclareMathOperator{\Aut}{Aut}

\DeclareMathOperator{\End}{End}

\DeclareMathOperator{\ram}{ram}
\DeclareMathOperator{\Gal}{Gal}

\newcommand{\PSL}[2]{\mathrm{PSL}({#1}, {#2})}
\newcommand{\PGL}[2]{\mathrm{PGL}({#1}, {#2})}
\newcommand{\SL}[2]{\mathrm{SL}({#1}, {#2})}
\newcommand{\GL}[2]{\mathrm{GL}({#1}, {#2})}
\newcommand{\GLtwo}{\mathrm{GL}(2)}
\newcommand{\abs}[1]{{\left| #1 \right|}}
\newcommand{\norm}[1]{{\|  #1 \|}}

\newcommand{\spandef}[2]{{  \left\langle  {#1}  \ \left| \   {#2} \right. \right\rangle }}
\newcommand{\setdef}[2]{{  \left\{  {#1}  \ \left| \   {#2} \right. \right\} }}
\newcommand{\setdefm}[3]{{  #1\{  {#2}  \ #1\vert \   {#3} #1\} }}
\newcommand{\mtext}[1]{\quad\mbox{#1}\quad}

\newcommand{\MCG}{\mathcal{MC\kern0.04emG}}
\newcommand{\ML}{\mathcal{M \kern0.07emL}}
\newcommand{\PML}{\mathcal{P \kern0.07em M \kern0.07em L}}

\renewcommand{\tilde}{\widetilde}
\newcommand{\num}[1]{\abs{#1}}

\newcommand{\h}{\mathcal H}
\newcommand{\A}{\mathbb A}

\newcommand{\oO}{\mathcal O}
\newcommand{\sQ}{\mathcal Q}
\newcommand{\p}{P}
\newcommand{\pbar}{\overline{P}}
\newcommand{\q}{Q}
\newcommand{\qbar}{{\overline{Q}}}
\newcommand{\qbarsquare}{{\qbar\vphantom{\q}^2}}

\newcommand{\hthree}{{\mathcal{H}^3}}
\newcommand{\sP}{{\mathcal P}}
\newcommand{\Qbar}{{\overline{\mathbb{Q}}}}
\newcommand{\FF}{\mathbb{F}}
\newcommand{\coleq}{\mathrel{\mathop:}=}
\newcommand{\newforms}{^{\mathrm{new}}}
\newcommand{\oldforms}{^{\mathrm{old}}}
\newcommand{\doublequom}[5]{{\raisebox{-#4}{$#1$}}#5\backslash {\raisebox{#4}{$#2$}}#5 / {\raisebox{-#4}{$#3$}}}

\newcommand{\rightquom}[4]{{\raisebox{#3}{$#1$}}#4/ {\raisebox{-#3}{$#2$}}}
\newcommand{\leftquom}[4]{{\raisebox{-#3}{$#1$}}#4\backslash{\raisebox{#3}{$#2$}}}

\newcommand{\mysmallmatrix}[4]{ \bigl( \begin{smallmatrix}  #1&#2\\ #3&#4 \end{smallmatrix} \bigr)
}

\newcommand{\Nrd}{\mathit{Nrd}}

\newcommand{\Norm}{{\mathcal{N}}}
\newcommand{\disc}{\mathit{disc}}

\newcommand{\cont}{\mathrm{cont}}

\newcommand{\hyp}{\nobreakdash-\hspace{0pt}}
\newcommand{\3}[1]{3\hyp}

\makeatletter
\DeclareRobustCommand*\tss[1]{%
  \@textsubscript{\selectfont#1}}
\def\@textsubscript#1{%
  {\m@th\ensuremath{_{\mbox{\fontsize\sf@size\z@#1}}}}}
\makeatother

\makeatletter
\@ifundefined{swappedhead@plain}{%
\def\swappedhead#1#2#3{%
   \thmnumber{\@upn{\@secnumfont#2\@ifnotempty{#1}{.~}}}%
   \thmname{#1}%
    \thmnote{ {\the\thm@notefont(#3)}}}
\let\swappedhead@plain=\swappedhead}
\makeatother

\theoremstyle{plain}

\swapnumbers
\newtheorem{theorem}{Theorem}[section]

\newtheorem{lemma}[theorem]{Lemma}

\newtheorem{proposition}[theorem]{Proposition}
\newtheorem{claim}[theorem]{Claim}
\newtheorem{question}[theorem]{Question}

\newtheorem{vfconj}[theorem]{Virtual Fibration Conjecture}
\newtheorem{forwardlemma}{Lemma}

\theoremstyle{definition}

\theoremstyle{remark}
\newtheorem{remark}[theorem]{Remark}

\makeatletter
  \let\c@theorem=\c@subsection
  \let\c@figure=\c@subsection
  \let\p@figure=\p@subsection
  
  \let\cl@figure=\cl@subsection
  \let\c@table=\c@subsection
  \let\p@table=\p@subsection
  
  \let\cl@table=\cl@subsection
  \let\c@equation=\c@subsection
  \let\p@equation=\p@subsection
  
  \let\cl@equation=\cl@subsection
\makeatother

\newenvironment{xyoverpic*}[3]
{%
\begin{xy}
\xyimport#1{\includegraphics[#2]{#3}}
}{\end{xy}}

\begin{document}

\title[Increasing the number of fibered faces of hyperbolic 3-manifolds]
{Increasing the number of fibered faces
 of \\  arithmetic hyperbolic 3-manifolds}

\author[Dunfield]{Nathan M.~Dunfield}
\address{Department  of Mathematics, MC-382, University of Illinois, Urbana, IL 61801, USA}
\email{nmd@uiuc.edu}

\thanks{Both authors were partially supported by the NSF, and Dunfield
  was partially supported by the Sloan Foundation. \hfill  \versioninfo}

\author[Ramakrishnan]{Dinakar Ramakrishnan} \address{Mathematics 253-37, California Institute of Technology, Pasadena, CA 91125, USA}

\email{dinakar@caltech.edu}

\begin{abstract}
  We exhibit a closed hyperbolic 3-manifold which satisfies a very
  strong form of Thurston's Virtual Fibration Conjecture.  In
  particular, this manifold has finite covers which fiber over the
  circle in arbitrarily many ways.   More precisely, it has a
  tower of finite covers where the number of fibered faces of the
  Thurston norm ball goes to infinity, in fact faster than any power
  of the logarithm of the degree of the cover, and we give a more
  precise quantitative lower bound.  The example manifold $M$ is
  arithmetic, and the proof uses detailed number-theoretic
  information, at the level of the Hecke eigenvalues, to drive a
  geometric argument based on Fried's dynamical characterization of
  the fibered faces. The origin of the basic fibration $M\to S^1$ is
  the modular elliptic curve $E=X_0(49)$, which admits multiplication
  by the ring of integers of $\Q[\sqrt{-7}]$. We first base change the
  holomorphic differential on $E$ to a cusp form on $\GLtwo$ over
  $K=\Q[\sqrt{-3}]$, and then transfer over to a quaternion algebra
  $D/K$ ramified only at the primes above $7$; the fundamental group
  of $M$ is a quotient of the principal congruence subgroup of
  $\oO_D^\ast$ of level $7$.  To analyze the topological properties of
  $M$, we use a new practical method for computing the Thurston norm,
  which is of independent interest.  We also give a non-compact
  finite-volume hyperbolic 3-manifold with the same properties by
  using a direct topological argument.
 \end{abstract}

%

\maketitle

\setcounter{tocdepth}{1}
\tableofcontents

\thispagestyle{empty}



\section{Introduction}
\label{sec_intro}

The most mysterious variant of the circle of questions surrounding
Waldhausen's Virtual Haken Conjecture \cite{Waldhausen68} is:
\begin{vfconj}[Thurston]\label{conj_virt_fib}
  If $M$ is a finite-volume hyperbolic \3-manifold, then $M$ has a
  finite cover which fibers over the circle, i.e.~is a surface bundle
  over the circle.
\end{vfconj}
This is a very natural question, equivalent to asking whether
$\pi_1(M)$ contains a \emph{geometrically infinite} surface group.
However, compared to the other forms of the Virtual Haken
Conjecture, there are relatively few non-trivial examples where it
is known to hold, especially in the case of closed manifolds (but
see \cite{Reid:bundles, Leininger2002, Walsh2005, Button2005} and
especially \cite{Agol2007}).  Moreover, there are indications that
fibering over the circle is, in suitable senses, a rare property
compared, for example, to simply having non-trivial first cohomology
\cite{DunfieldDThurston, Masters2005}.

Despite this, we show here that certain manifolds satisfy
Conjecture~\ref{conj_virt_fib} in a very strong way, in that they have
finite covers which fiber over the circle in many distinct ways.  For
a \3-manifold $M$, the set of classes in $H^1(M; \Z)$ which can be
represented by fibrations over the circle are organized by the
Thurston norm on $H^1(M; \R)$.  The unit ball in this norm is a finite
polytope where certain top-dimensional faces, called the fibered
faces, correspond to those cohomology classes coming from
fibrations (see Section~ \ref{sec_fiberedfaces} for details).  The
number of fibered faces thus measures the number of fundamentally
different ways that $M$ can fiber over the circle.  We will sometimes
abusively refer to these faces of the Thurston norm ball as ``the
fibered faces of $M$''.

If $N \to  M$ is a finite covering map, the induced map $H^1(M; \R) \to
H^1(N; \R)$ takes each fibered face of $M$ to one in $N$; if $H^1(N; \R)$ is
strictly larger than $H^1(M; \R)$, then it may (but need not) have
additional fibered faces.   The qualitative form of our main result is:
\begin{theorem}\label{intro_thm}
  There exists a closed hyperbolic \3-manifold $M$ which has a
  sequence of finite covers $M_n$ so that the number of fibered faces
  of the Thurston norm ball of $M_n$ goes to infinity.
\end{theorem}
Moreover, we prove a quantitative refinement of this result
(Theorem~\ref{rapid_growth}) which bounds from below the number of
fibered faces of $M_n$ in terms of the degree of the cover. While it
is the closed case of Conjecture~\ref{conj_virt_fib} that is most
interesting, we also give an example of a \emph{non-compact}
finite-volume hyperbolic 3-manifold with the same property
(Theorem~\ref{thm_whitehead}) using a simple topological
argument.

The example manifold $M$ of Theorem~\ref{intro_thm} is arithmetic,
and the proof uses detailed number-theoretic information about it,
at the level of the Hecke eigenvalues, to drive a geometric argument
based on Fried's dynamical characterization of the fibered faces. To
state the geometric part of the theorem, we need to introduce the
Hecke operators (see Section~\ref{sec_hecke_ops} for details).
Suppose $M$ is a closed hyperbolic \3-manifold, and we have a pair
of finite covering maps $p, q \maps N \to M$; when $M$ is arithmetic
there are many such pairs of covering maps because the commensurator
of $\pi_1(M)$ in $\Isom(\hthree)$ is very large.  The associated
\emph{Hecke operator} is the endomorphism of $H^1(M)$ defined by
$T_{p,q} = q_* \circ p^*$, where $q_* \maps H^1(N) \to H^1(M)$ is
the transfer map.  The simplest form of our main geometric lemma is
the following:
\begin{lemma}\label{basic_geom_thm}
  Let $M$ be a closed hyperbolic 3-manifold, and $p, q \maps N \to M$ a
  pair of finite covering maps.  If $T_{p,q}(\omega) = 0$ for some $\omega \in
  H^1(M;\Z)$ coming from a fibration over the circle, then $p^*(\omega)$
  and $q^*(\omega)$ lie in distinct fibered faces.
\end{lemma}

We prove this lemma in Section~\ref{sec-geom-thm}. Then in
Section~\ref{sec_example}, we give a manifold with an infinite tower
of covers to which Lemma~\ref{basic_geom_thm} applies at each step,
thus proving Theorem~\ref{intro_thm}.  In fact, we show the
following refined quantitative version:
\begin{theorem}\label{rapid_growth}
  There is a closed hyperbolic 3-manifold $M$ of arithmetic type, with
  an infinite family of finite covers $\{M_n\}$ of degree $d_n$, where
  the number $\nu_n$ of fibered faces of $M_n$ satisfies
  \[
  \nu_n \, \geq \, \exp\left(0.3\frac{\log d_n}{\log\log d_n}\right) \quad \mbox{as $d_n \to \infty$.}
  \]
  In particular, for any $t <1$, there is a constant $c_t>0$ such
  that
  \[
  \nu_n \, \geq \, c_t e^{(\log d_n)^t}.
  \]
\end{theorem}
Note that this bound for $\nu_n$ is
slower than any positive power of $d_n$, but is faster than any
(positive) power of $\log d_n$. For context, the Betti number
$b_1(M_n) = \dim H^1(M_n; \R)$ is bounded above by (a constant
times) the degree $d_n$, and bounded below (see Proposition
\ref{prop_lower_bound}), for any $\varepsilon
>0$, by (a constant times) $d_n^{1/2-\varepsilon}$ for $n$ large
enough (relative to $\varepsilon$), while $\nu_n$ is at least as large as (a constant
times) $e^{(\log d_n)^{0.99}}$. In our non-compact example of
Theorem~\ref{thm_whitehead}, the number of faces grows exponentially
in the degree while $b_1$ grows linearly.

We now describe the basic ideas behind the construction of the
manifolds in Theorem~\ref{rapid_growth}.  For an arithmetic hyperbolic
\3-manifold $M$, one has a Hecke operator as above associated to each
prime ideal of the field of definition.  The key to applying
Lemma~\ref{basic_geom_thm} repeatedly is to have a fibered class $\omega \in
H^1(M)$ which is killed by infinitely many such Hecke operators.  One
can produce cohomology classes which are killed by infinitely many
Hecke operators using the special class of CM forms.  Fortunately,
there is a manifold of manageable size whose cohomology contains a CM
form coming from base change of an automorphic form associated to a
certain elliptic curve with complex multiplication, and that class
turns out to fiber over the circle!

\subsection{Outline  of the arithmetic construction}
\label{sketch_of_construction}

Here is a sketch of how the manifolds $\{M_n\}$ of
Theorem~\ref{rapid_growth} are built from arithmetic; for details, see
Section~\ref{sec_example}.     Let $E$ be the elliptic curve over
$\Q$ defined by $y^2+xy=x^3-x^2-2x-1$, which has conductor $49$ and
admits complex multiplication by the ring of integers of the
imaginary quadratic field $\Q[\sqrt{-7}]$. It corresponds to a
holomorphic cusp form of weight $2$ for the congruence subgroup
$\Gamma_0(49)$ of $\SL{2}{\Z}$ acting on the upper half plane
${\mathcal H}$, given by $f(z) = \sum_{n\geq 1} a_nq^n$ with
$q=e^{2\pi iz}$, where for every prime $p \neq 7$, the eigenvalue of
$f$ under the Hecke operator $T_p$ is $a_p=p+1-\num{E({\mathbb
F}_p)}$. The differential $f(z)dz$ is invariant under
$\Gamma_0(49)$, and hence defines a holomorphic $1$-form on
$Y_0(49):=\Gamma_0(49)\backslash{\mathcal H}$, and by the
cuspidality of $f$, this differential extends to the natural (cusp)
compactification of $Y_0(49)$.
Put $K=\Q[\sqrt{-3}]$, in which the prime $7$ splits as $\q
{\qbar}$, with $\q =\left(2+\big(1+i\sqrt{3}\big)\big/2\right)$. Let
$f_K$ denote the base change of $f$ to $K$, which is a cusp form on
the hyperbolic $3$-space $\hthree$ of ``weight $2$'' and level
$\q^2\qbarsquare$
for the group
$\GL{2}{\oO_K}$. One can associate to $f_K$ a cuspidal automorphic
representation $\pi'$ of $\GL{2}{\A_K}$ with trivial central
character and conductor $\q^2\qbarsquare$; here $\A_K$ denotes the
ad\`ele ring of $K$. By a basic property of base change, we have,
for every degree $1$ prime $P$ of $\oO_K$ above a rational prime
$p\neq 7$ unramified in $K$, the $P^{\mathrm{th}}$ Hecke eigenvalue $a_P(f_K)$
of $f_K$ equals the $p^{\mathrm{th}}$ Hecke eigenvalue of $f$, namely
$a_p=p+1-\num{E({\mathbb F}_p)}$.

Let $D$ be the quaternion division algebra over $K$ which is
ramified exactly at the primes $\q$ and $\qbar$. Fortunately for us,
the local components of $\pi'$ at $\q$ and $\qbar$ are both
supercuspidal (see Lemma 5.3); this has to do with the fact that $E$
does not acquire good reduction over an abelian extension of $\Q_7$,
this fact being controlled, thanks to a useful criterion of
D.~Rohrlich \cite{Rohr}, by the valuation (at $\q$ and $\qbar$) of
the discriminant $\Delta$ of $E$. (After writing this paper, we
learned of an earlier proof of the supercuspidality at $7$ in
\cite{GL}, where the argument is somewhat different.) Therefore, by
the Jacquet-Langlands correspondence, there is an associated cusp
form $h$ of weight $2$ on $\hthree$ relative to a congruence
subgroup $\Gamma$ of the units in a maximal order $\oO_D$, such that
for every unramified prime $P\neq \q,\qbar$, the $P^{\mathrm{th}}$
Hecke eigenvalues of $f_K$ and of $h$ coincide. Moreover, for $P\in
\{\q, \qbar\}$, we can read off the conductor and the dimension of
the associated representation of $(D\otimes_K K_P)^\ast$ from the
local correspondence (see Lemma 5.5). It follows that $\Gamma$,
which is co-compact, is the principal congruence subgroup of level
$7= \q \qbar$. Our base manifold in Theorem~\ref{rapid_growth} is
$M=X(7)=\Gamma\backslash \hthree$. While $H^1(M)$ is \3-dimensional,
we show that the new subspace of $H^1(M)$ is 2-dimensional and, as a
module under the Hecke algebra of correspondences, is isotypic of
type $\eta=\eta_h$, the cohomology class defined by $h$. A difficult
computation shows that $M$ fibers over the circle, and moreover,
this can be associated with a class of type $\eta$.

Now consider the set $\sP$ of rational primes $p\neq 7$ which are inert
in $\Q[\sqrt{-7}]$, but are split in $K$ as $\p \pbar$; it has density
$1/4$. Then, if we put
\[
\sP_K= \setdef{P}{\Norm_{K/\Q}(P)=p\in \sP},
\]
the Hecke operators $T_P$ act by zero on $\eta$. Each $P$ in $\sP_K$
gives a covering $M(P)$ of $M$ of degree $\Norm(P) +1$ which defines
the associated Hecke operator.  Lemma~\ref{basic_geom_thm} allows us, since
$T_P\eta=0$, to conclude that the two natural transforms $\eta_1, \eta_P$ of
$\eta$ define cohomology classes of $M(P)$ which lie on two different
fibered faces.  If we order $\sP_K$ according to the rational primes
$p$ defining them, then we inductively build a covering
$M_n=M(P_1\dots P_n)$ of degree $d_n\coleq\prod_{j=1}^n(1+p_j)$ such that
there are at least $2^n$ distinct fibered faces on $M_n$.  Using the
density of $\sP$, we get the lower bound for $\nu_n \geq 2^n$ given in
Theorem~\ref{rapid_growth} in terms of $d_n$.

To understand why we chose to look at this particular example, it
may be helpful to note the following. Suppose $D_0$ is an indefinite
quaternion algebra over $\Q$, and $\Delta$ a congruence subgroup of
$D_0^\ast$ with associated \emph{Shimura curve} $S=S_\Delta$ over
$\C$, which is a compact Riemann surface. For any imaginary
quadratic field $K$ such that $D=D_0\otimes_\Q K$ is still a
division algebra, we may consider the hyperbolic $3$-manifold
$M=\Gamma\backslash\h_3$, for a congruence subgroup $\Gamma$ of
$D^\ast$. When $\Gamma \cap D_0^\ast = \Delta$, the surface $S$
embeds in $M$, and since it is totally geodesic, the cohomology
class defined by $S$ cannot give rise to any fibering of $M$ over
the circle \cite{ThurstonFibered}. This suggests that we will fail
to construct cohomology classes on $M$ with the desired fibering
property if we transfer to $D^\ast$ a cusp form on $\GLtwo/K$ which
is the base change $f_K$ of an elliptic modular cusp form $f$ of
weight $2$ which will transfer to such a $D_0^\ast$. It is not an
accident that we chose our example above where the $f$ of interest
does not transfer to any indefinite quaternion algebra over $\Q$.

Finally, it should be noted that if one starts with a non-CM
elliptic curve $E$ over $\Q$, then one knows, by a theorem of
Elkies, that there are infinitely many primes $p$ for which $a_p$ is
zero \cite{Elkies1987}. But for our method to work we would also
need an example where this property holds for an infinite set $\sP$
of primes $p$ which split in a suitable imaginary quadratic field
$K$. Even then it would only give a qualitative result as $\sP$
would have density zero. Our quantitative result
(Theorem~\ref{rapid_growth}) depends on the density of the
corresponding $\sP$ in the CM case being $1/4$.

\subsection{Moral}

For an \emph{arithmetic} hyperbolic 3-manifold, the commensurator of
its fundamental group is very large, in fact dense in
$\Isom^+(\hthree) \cong \PSL{2}{\C}$.  Recently, there has been much
important work which exploits this density \emph{geometrically}, see
\cite{LackenbyLongReid2006, CooperLongReid2006, Venkataramana2006,
  Agol2006}.  On the number-theoretic side, the theory of automorphic
forms tells us a great deal about the cohomology of arithmetic
hyperbolic 3-manifolds as it relates to virtual Haken type questions
(see e.g.~\cite{Clozel1987, DunfieldFCalegari2006}).  To the best of
our knowledge, our work here is the first time that \emph{more
  refined} automorphic information has been combined with geometric
arguments and yields, for example, a geometric/topological
interpretation of the vanishing of the Hecke eigenvalues. Thus we
hope for deeper connections between these two areas in the future.
In particular, it would be very interesting to answer the following:
\begin{question}
  Is there an automorphic criterion which implies that certain
  cohomology classes of arithmetic hyperbolic 3-manifold give
  fibrations over the circle?
\end{question}

\subsection{Practical methods for computing the Thurston norm}

The example manifold of Theorem~\ref{rapid_growth} is quite
complicated from a topological point of view; its hyperbolic volume
is about 100 and triangulations of it need some 130 tetrahedra.
Despite this, we were able to compute its Thurston norm and check
that it fibers over the circle, which is necessary for the proof of
Theorem~\ref{rapid_growth}.  To do this, we used new methods for
both these tasks.  While loosely based on normal surfaces, these
techniques eschew guaranteed termination in favor of quick results.
The basic idea is to consider only normal surfaces representing
elements of $H^1(M)$ that are ``obvious'' with respect to the
triangulation, and then randomize the triangulation until a minimal
norm surface is found.  These same techniques have been useful in
many other examples and are of independent interest.  See
Sections~\ref{subsec_ex_thurston} and \ref{subsec-fibering-M} for a
complete description of our method, which can often determine
whether a manifold made up of several hundred tetrahedra fibers.
Computing the Thurston norm is more subtle, and in the case of our
particular $M$, we had to heavily exploit its symmetries, but we
also suggest a general method, as of yet untested, for attacking
this.

\subsection{Improvements}

In subsequent work, Long and Reid have considerably strengthened
Lemma~\ref{basic_geom_thm} and its extension Theorem~\ref{cong_cover_thm} by
showing that the hypothesis on the Hecke operators can be dispensed
with:
\begin{theorem}[Long and Reid \cite{LongReid2007}] \label{thm_long_reid}
  Let $M$ be a closed arithmetic hyperbolic 3-manifold.  If $M$ fibers
  over the circle, then $M$ has finite covers whose Thurston norm balls have
  arbitrarily many fibered faces.
\end{theorem}
In addition to the work of Fried \cite{Fried1979} on which
Lemma~\ref{basic_geom_thm} hinges, the proof of Theorem~\ref{thm_long_reid}
uses work of Cooper, Long, and Reid on suspension pseudo-Anosov flows
\cite{CooperLongReid1994}, as strengthened by Masters
\cite{Masters2006}.  In Section~\ref{sec_long_reid} we give a
simplified proof of Theorem~\ref{thm_long_reid} using only Fried's
theorem; a different, but equally concise, simplification was given by
Agol~\cite{Agol2007}.  In any case, the soft nature of its proof mean
that Theorem~\ref{thm_long_reid} cannot be used to prove quantitative
results such as Theorem~\ref{rapid_growth}.  However,
Theorem~\ref{thm_long_reid} does have the advantage that it is much
easier to apply.

In a major breakthrough, Agol has just shown there are infinitely many
commensurability classes of arithmetic hyperbolic 3-manifolds which
fiber over the circle \cite{Agol2007}.  Combining with
Theorem~\ref{thm_long_reid}, this means the qualitative behavior of
Theorem~\ref{intro_thm} actually occurs in an infinite number of
examples, providing further evidence for
Conjecture~\ref{conj_virt_fib}.

In Theorems~\ref{rapid_growth} and \ref{thm_long_reid}, the distinct
fibered faces are all equivalent under the isometry group of
the cover manifolds; indeed this is intrinsic to the method.  A natural
question is whether one can find a tower of covers where the fibered
faces fall into arbitrarily many classes modulo isometries.  In the
case of manifolds with cusps, Theorem~\ref{thm_whitehead} gives such
examples since the number of faces grows exponentially in the degree
of the cover, whereas the size of the isometry group of a hyperbolic
manifold is bounded linearly in the volume.  

\subsection{Paper outline}

In Section~\ref{sec_fiberedfaces}, we review the basics of the
Thurston norm and Fried's dynamical characterization of the fibered
faces.  In Section~\ref{sec-geom-thm}, we discuss Hecke operators and
congruence covers in the context of arithmetic hyperbolic
\3-manifolds, and then prove Lemma~\ref{basic_geom_thm} and
its generalization Theorem~\ref{cong_cover_thm} which underpins
Theorem~\ref{rapid_growth}.  We give the precise description of the
manifold used in Theorem~\ref{rapid_growth} in
Section~\ref{sec_example}.  We then analyze the automorphic and
topological properties of this manifold in
Sections~\ref{sec-automorphic} and \ref{sec-topological-props}
respectively.  In Section~\ref{sec_proof_of_main_thm} we assemble the
pieces and prove Theorem~\ref{rapid_growth}.
Section~\ref{sec-whitehead} gives our example of this phenomenon in
the case of hyperbolic 3-manifolds with cusps.  Finally,
Section~\ref{sec_long_reid} gives our simplified proof of
Theorem~\ref{thm_long_reid}.

\subsection{Acknowledgments}

Dunfield was partially supported by a Sloan Fellowship and US NSF
grants DMS-0405491 and DMS-0707136, and some of this work was done
while he was at Caltech.  Ramakrishnan was partially supported by US
NSF grants DMS-0402044 and DMS-0701089.  We thank Philippe Michel
for a useful discussion concerning the quantitative version of our
main result (Theorem~\ref{rapid_growth}), and Danny Calegari for
suggesting that we look at families of link complements in the
cusped case, which led us to Theorem~\ref{thm_whitehead}.  We also
thank Darren Long and Alan Reid for sending us an early version of
\cite{LongReid2007}, and Benedict Gross for informing us of an earlier
proof of Lemma 5.3 in \cite{GL}. Finally, we thank the authors of
the software packages \cite{SnapPea, Snap, Cremona1984,
  MAGMA214, SAGE}, without which this paper would not have been possible,
  as well as the referee
  for carefully reading the paper and suggesting improvements in exposition.

\section{Fibered faces of the Thurston norm ball}
\label{sec_fiberedfaces}

Let $M$ be a closed orientable hyperbolic 3-manifold.  When a
cohomology class $\omega \in H^1(M; \Z) ={\rm Hom}(\pi_1(M),\Z)$ can be
represented by a fibration $M \to S^1$, we say that it \emph{fibers}.
In this section, we review the work of Thurston and Fried on the
structure of the set of fibered classes in $H^1(M; \Z)$.  It is not
hard to see that $\omega \in H^1(M; \Z)$ fibers if and only if it can be
represented by a nowhere vanishing 1-form: a fibration gives rise to
such a form by pulling back the standard orientation 1-form on $S^1$,
and conversely such a form can be integrated (since its periods are
integers) into a map to $S^1$ which is a fibration.  If
we pass to real coefficients, the latter condition clearly defines
an \emph{open} subset $U$ of $H^1(M; \R) \setminus \{0\}$; a
fundamental result of Thurston shows that $U$ must have the
following very restricted form.  This set is prescribed by the
Thurston norm, which measures the simplest surface that represents
the Poincar\'e dual of a cohomology class.

More precisely, for $\omega \in H^1(M; \Z)$ define its \emph{Thurston norm}
by
\[
\| \omega \| = \min \setdefm{\big}{-\chi(\Sigma)}{\mbox{$\Sigma$ is an embedded orientable
    surface dual to $\omega$} },
\]
where $\chi(\Sigma)$ is the Euler characteristic of $\Sigma$, and we
further require that $\Sigma$ has no 2-sphere components.  Thurston
showed that this gives a norm on $H^1(M; \Z)$ which extends
continuously to one on $H^1(M; \R)$; see \cite{ThurstonNorm} for
details of the assertions made in this paragraph.
Moreover, the unit ball $B$ of this norm is a
bounded convex polytope, i.e.~the convex hull of finitely many
points. Moreover, there are top-dimensional faces of $B$, called the
\emph{fibered faces}, so that $\omega \in H^1(M;\Z)$ fibers if and
only if it lies in the cone over the interior of a fibered face,
that is, the ray from the origin through $\omega$ intersects the
interior of such a face.  We will say that such a fibered $\omega$
\emph{lies in} the corresponding fibered face. Finally, as the map
$\omega \mapsto -\omega$ preserves fibering, the fibered faces come
in pairs interchanged by this symmetry of the Thurston norm ball.
When determining the number of fibered faces, we often count in
terms of these fibered face pairs, and say that fibered classes
$\alpha$ and $\beta$ lie in \emph{genuinely distinct} fibered faces
if both $\alpha$ and $-\alpha$ are not in the fibered face of
$\beta$.

\subsection{Behavior of fibered faces under covers}\label{sec_norm_cover}

Now suppose $p \maps N \to M$ is a finite covering map.  The natural
map $p^* \maps H^1(M; \R) \to H^1(N; \R)$ is an embedding, and a
deep theorem of Gabai shows that $p^*$ preserves the Thurston norm
\cite[Cor.~6.18]{Gabai83}.  Equivalently, if we denote the unit
Thurston norm balls of $M$ and $N$ by
$B_M$ and $B_N$ respectively, we have $p^*(B_M) = B_N \cap
p^*\left(H^1(M; \R)\right)$.  By work of Stallings, a class $\omega
\in H^1(M; \Z)$ represents a fibration if and only if $p^*(\omega)$
does \cite{Stallings62}.  In particular, if $N$ fibers but $M$ does
not, then $H^1(N)$ is larger than $H^1(M)$.  Thus each fibered face
of $B_M$ gives rise to a distinct fibered face of $B_N$; if $N$ has
additional cohomology, we can hope that $B_N$ has new fibered faces,
but the interiors of these must be disjoint from
$p^*\left(H^1(M)\right)$.

\subsection{Fried's work}\label{subsec_frieds_work}

We now turn to Fried's dynamical characterization of when two
fibrations lie in the same fibered face, which is in terms of a
certain flow that is transverse to the fibers of the fibration $M
\to S^1$.  Suppose $\phi$ is a self-homeomorphism
of a closed surface $\Sigma$ of genus
at least 2, and consider the mapping torus with monodromy $\phi$:
\[
M_\phi = \raisebox{4pt}{$\Sigma \times [0,1]$} \Big/ \raisebox{-4pt}{$(s,1) \sim
 (\phi(s),0)$}
\]
Thurston proved that $M_\phi$ is hyperbolic if and only if $\phi$ is
what is called pseudo-Anosov \cite{ThurstonFibered, Otal96}.  The
latter means that $\phi$ is isotopic to a homeomorphism which
preserves a pair of foliations of $\Sigma$ in a controlled way.
Henceforth, we always assume that $\phi$ has been isotoped
to such a preferred representative.
Now, $M_\phi$ has a natural suspension flow $\F_\phi$ which is
transverse to the circle fibers, where a point moves at unit speed
in the $[0,1]$--direction. We will call $\F_\phi$ the
\emph{transverse pseudo-Anosov flow}.

Conversely, given $\omega \in H^1(M; \Z)$ coming from a fibration, the
monodromy $\phi$ of the bundle structure is well-defined up to isotopy.
Thus there is a corresponding transverse pseudo-Anosov flow $\F_\omega$,
which is well-defined up to isotopy.  Fried's first result is:
\begin{theorem}[Fried \cite{Fried1979}]\label{thm_fried}
  Let $M$ be a closed orientable hyperbolic \3-manifold.  Then two
  fibrations of $M$ over the circle lie over the same fibered face if and
  only if the corresponding transverse pseudo-Anosov flows are isotopic.
\end{theorem}
Fried also provided the following characterization of those $\omega$ lying
over a fibered face.  Let $\F$ be the flow associated to a
particular fibered face $F$, and let $D \subset H_1(M; \R)$ be the set of
\emph{homology directions} of $\F$, namely the set of all accumulation
points of the homology classes of long, nearly closed orbits of $\F$.
Fried showed that the dual cone to $D$,
\[
C = \setdef{\omega \in H^1(M; \R)}{ \mbox{$\omega(v) > 0$ for all $v \in D$}},
\]
is precisely the cone on the interior of the fibered face $F$
\cite[Thm.~7]{Fried1979}.

One kind of element of $D$ is the homology class of a
closed flowline of $\F$.  There are always closed flowlines, for
instance coming from the orbits of the finitely many singular points of
the invariant foliations, which are permuted
among themselves by the monodromy.  The particular consequence of Fried's
work that we will use here, which is immediate from our discussion, is:
\begin{lemma}[Fried]\label{disjointness_crit}
  Let $M$ be a closed orientable hyperbolic \3-manifold, with fibered
  classes $\alpha, \beta \in H^1(M; \Z)$.  Let $c$ be a closed orbit of the
  transverse pseudo-Anosov flow for $\alpha$.  If $\beta(c) = 0$, then $\alpha$
  and $\beta$ lie over genuinely distinct fibered faces.
\end{lemma}

\section{The geometric theorem}\label{sec-geom-thm}

We begin this section by proving Lemma~\ref{basic_geom_thm}, which
contains the key geometric idea of this paper: a fibered cohomology
class which is annihilated by a Hecke operator gives rise to two
genuinely distinct fibered faces in the corresponding cover.  We then
give a detailed review of Hecke operators for arithmetic hyperbolic
3-manifolds, and end by proving the complete version of our geometric
result (Theorem~\ref{cong_cover_thm}) which is needed to prove
Theorem~\ref{rapid_growth}.

\subsection{Hecke operators}
\label{sec_hecke_ops}

From a geometric point of view, a Hecke operator is the map on
cohomology induced from the following setup.  Suppose $M$ is a
topological space, and we have a pair of \emph{finite} covering maps
$p_, q \maps N \to M$.  The map on singular chains $C_*(M) \to C_*(N)$
which takes a singular simplex to the sum of its inverse images under
$q$ induces transfer homomorphisms $H_*(M) \to H_*(N)$ and $H^*(N) \to
H^*(M)$ which run in opposite directions to the usual maps that $q$
induces on (co)homology (see e.g.~\cite[\S3.G]{HatcherBook} for
details).  The Hecke operator $T_{p,q}$ of this pair of covering maps
is the endomorphism of $H^*(M)$ defined by $q_* \circ p^*$, where $p^*
\maps H^*(M) \to H^*(N)$ the pullback map and $q_* \maps H^*(N) \to
H^*(M)$ is the transfer map.
\begin{remark}
  In this paper, $M$ will always be a \3-manifold and we will be
  interested in the Hecke operator on $H^1(M)$.  While the argument
  below is given purely in terms of cohomology, the geometrically
  minded reader may prefer to contemplate the Poincar\'e dual
  group $H_2(M)$.  There, the Hecke operator $T_{p_, q}$ on $H_*(M)$
  is the composite $q_* \circ p^*$, where $p^*$ is the transfer map.
  Such Hecke operators commute with the Poincar\'e duality isomorphism
  $H^1(M) \cong H_2(M)$, and so it makes no difference whether one takes
  the homological or cohomological point of view.  The action of
  $T_{p,q}$ on a class $\omega \in H_2(M)$ is particularly simple to think of
  geometrically: If an embedded surface $\Sigma \subset M$ represents $\omega$,
  then the immersed surface $q(p^{-1}(\Sigma))$ represents $T_{p,q}(\omega)$.
\end{remark}

\subsection{Main geometric idea}
When $M$ is arithmetic, there are many manifolds $N$ which cover it in
distinct ways.  Before getting into this, let us give the central
topological idea of Theorem~\ref{intro_thm} in its simplest form from
the introduction:

\begin{forwardlemma}
  Let $M$ be a closed hyperbolic 3-manifold, and suppose $\omega \in H^1(M;
  \Z)$ comes from a fibration of $M$ over the circle.  Further assume
  that $p, q \maps N \to M$ are a pair of finite covering maps.  If
  $T_{p_, q}(\omega) = 0$ then $p^*(\omega)$ and $q^*(\omega)$ lie in genuinely
  distinct fibered faces.
\end{forwardlemma}

\begin{proof}
  We will apply Lemma~\ref{disjointness_crit} to justify our claim.
  Downstairs in $M$, let $c$ be a closed orbit of the pseudo-Anosov
  flow associated to the fibration $\omega$.  Then $q^{-1}(c)$ is a closed
  orbit of the flow associated to the fibration $q^*(\omega)$.  To calculate
  $p^*(\omega)(q^{-1}(c))$, note that for any $\alpha \in H^1(N)$ one has
  $\alpha(q^*(c)) = \left( q_*(\alpha) \right)(c)$ and hence
  \[
  p^*(\omega)\left(q^{-1}(c)\right) = p^*(\omega)\left(q^*(c)\right) =
  q_*\left(p^*(\omega)\right)(c) = T_{p,q}(\omega)(c) = 0
  \]
  as required to apply the lemma.
\end{proof}
In the case of a tower of covers of $M$, the following strengthening
of the previous lemma will be needed to work inductively:
\begin{lemma}\label{full_geom_thm}
  Let $M$ be a closed hyperbolic 3-manifold, and $p, q \maps N \to M$
  be a pair of finite covering maps.  Suppose $\omega_1, \omega_2, \ldots ,\omega_n \in
  H^1(M; \Z)$ lie in genuinely distinct fibered faces.  If $T_{p_,
    q}(\omega_i) = 0$ for all $i$, then $N$ has at least $2 n$ pairs of
  fibered faces.  More precisely, $\left\{ p^*(\omega_1), \ldots , p^*(\omega_n),
  q^*(\omega_1), \ldots , q^*(\omega_n) \right\}$ lie in genuinely distinct fibered faces.
\end{lemma}

\begin{proof}
  By the discussion in Section~\ref{sec_norm_cover}, it is clear that
  the $p^*(\omega_i)$ lie in distinct fibered faces, as do the
  $q^*(\omega_j)$.  Thus it remains to distinguish the face of $p^*(\omega_i)$
  from that of $q^*(\omega_j)$.  As before, let $c$ be a closed orbit of the
  pseudo-Anosov flow associated to $\omega_j$, so that $q^{-1}(c)$ is an
  orbit of the flow associated to $q^*(\omega_j)$.  Then
  \[
  p^*(\omega_i)\left(q^{-1}(c)\right) = p^*(\omega_i)\left(q^*(c)\right) = q_*\left(p^*(\omega_i)\right)(c) =  T_{p,q}(\omega_i)(c) = 0
  \]
  and so Lemma~\ref{disjointness_crit} applies as needed.
\end{proof}

\subsection{Sources of Hecke operators}\label{sec-source-hecke}

Now we turn to the source of such multiple covering maps.  For a
hyperbolic \3-manifold $M$, let $\Gamma$ be its fundamental group, thought
of as a lattice in $\PSL{2}{\C}$.  The commensurator of $\Gamma$ is the
subgroup
\[
\Comm(\Gamma) = \setdef{g \in \PSL{2}{\C}}{\mbox{$g^{-1} \Gamma g \cap \Gamma$ is
    finite index in both $\Gamma$ and $g^{-1} \Gamma g$}}.
\]
When $M$ is arithmetic, $\Comm(\Gamma)$ is dense in $\PSL{2}{\C}$, and
Margulis showed that the converse is true as well; indeed, if $M$ is not
arithmetic, then $\Gamma$ has finite index in $\Comm(\Gamma)$.

Regardless, $g \in \Comm(\Gamma)$ can be  associated to a Hecke
operator as follows.  Let $\Gamma_g = g^{-1} \Gamma g \cap \Gamma$,
and $M_g = \Gamma_g \backslash \hthree$ be the corresponding closed
hyperbolic \3-manifold. Consider the finite covering map $p_g \maps
M_g \to M$ induced by the inclusion $\Gamma_g \to \Gamma $. We will
define a second such covering map by considering $\Gamma_{g^{-1}} =
g \Gamma g^{-1} \cap \Gamma$, and analogously $p_{g^{-1}} \maps
M_{g^{-1}} \to M$.  Now, as $g \Gamma_{g} g^{-1} = \Gamma_{g^{-1}}$,
the action of $g$ on $\hthree$ induces an isometry $\tau_g \maps M_g
\to M_{g^{-1}}$ giving us the picture
\[
\xymatrix@C-=2ex{
  M_g  \ar[rr]^{\tau_g} \ar[dr]_{p_g} & &M_{g^{-1}} \ar[dl]^{p_{g^{-1}}} \\
   & M
}
\]
Combining, we get a covering map $q_g \maps M_g \to M$ by taking
$q_g = p_{g^{-1}} \circ \tau_g$.  On the level of groups, the
covering map $q_g$ corresponds to the homomorphism $\Gamma_g \to
\Gamma$ given by $\gamma \mapsto g\gamma g^{-1}$. By definition, the
Hecke operator associated to $g$ is $T_{p_g, q_g}$ in the notation
of Section~\ref{sec_hecke_ops}.

\subsection{Standard congruence covers}
\label{cong_covers}

In number theory, one usually restricts to a subset of all the Hecke
operators described above which have nice collective properties.  For
those readers who are unfamiliar with this, and to explain how it
interacts with Lemma~\ref{full_geom_thm}, we give a detailed review of the
basic setup in Sections~\ref{cong_covers}--\ref{sec_new_and_old}.
This material is standard, and those with a background in number
theory may wish to skip ahead to Section~\ref{sec_cong_cover_thm}.

To describe these, we first review the general arithmetic construction
of a lattice in $\Isom^+(\hthree)$; for details see
e.g.~\cite{VignerasBook,MaclachlanReidBook}.  Begin with a number
field $K$ with exactly one complex place, and choose a quaternion
algebra $D$ over $K$ which is ramified at all real places of $K$.  Now
at the complex place of $K$, the algebra $D \otimes_K \C$ must be
isomorphic to $M_2(\C)$, the algebra of $2 \times 2$ matrices.  Taking the
units gives a homomorphism $D^\ast \hookrightarrow \GL{2}{\C}$; dividing out by the
respective centers embeds $D^\ast/K^\ast$ as a dense subgroup of
$\PGL{2}{\C} \cong \PSL{2}{\C} \cong \Isom^+(\hthree)$. Now if $\oO_D$ is a
maximal order of $D$, then the image $\Gamma$ of the units $\oO_D^\ast$ in
$\PGL{2}{\C}$ is a lattice, which is cocompact provided $D$ is a
division algebra rather than $M_2(K)$. We will denote the quotient
hyperbolic \3-orbifold as $X = \Gamma \backslash \hthree$.

For each ideal $A$ of $\oO_K$, we can define a corresponding
congruence orbifold $X(A)$ covering $X$ as follows.  Since we will
only need that case, we restrict to when $A$ is square-free.  The
key case is that of a prime ideal $\p$. Let $K_\p$ be the local
completion of $K$ at the
place $\p$; its valuation ring is denoted $\oO_{\p}$ with maximal
ideal $\sP$ and residue field $\FF = \oO_\p/\sP$.  The cover $X(\p)$
is constructed using the local algebra $D_\p = D \otimes_K K_\p$.

First, suppose $D$ does not ramify at $P$, i.e.~$D_\p \cong M_2(K_\p)$.
The cover $X(\p)$ is the congruence cover of ``$\Gamma_0$-type'' built as
follows.  The order $\oO_{D_P} \coleq \oO_D \otimes_{\oO_K} \oO_P$ is maximal
in $D_\p$, and so is conjugate to $M_2(\oO_\p)$.  Thus we get
\[
\oO^\ast_D \hookrightarrow \oO_{D_\p}^\ast \cong \GL{2}{\oO_p} \to
\rightquom{\oO_{D_\p}^\ast}{1 + \sP \oO_{D_P}}{3pt}{\big} \cong \GL{2}{\FF}
\]
which induces a homomorphism $\Gamma \to \PGL{2}{\FF}$.  By strong
approximation, the image of $\Gamma$ acts transitively on
$P^1(\FF)$.  By definition, $X(\p) = \Gamma_0(\p) \backslash
\hthree$ where $\Gamma_0(\p)$ is the $\Gamma$-stabilizer of a point
in $P^1(\FF)$.  This gives a cover $X(\p) \to X$ of degree
$\num{P^1(\FF)} = N(\p) + 1$, where $N(\p)=\vert\FF\vert$ is the
norm of $\p$.

Suppose instead $\p$ is one of the finitely many primes where $D$
ramifies.  Then $D_\p$ is the unique quaternion division algebra over
$K_\p$, and so is a local skew-field whose valuation ring is
$\oO_{D_\p}$ with a unique maximal bi-ideal $\sQ$.  If $q =
\num{\FF}$ we have
\[
\oO^\ast_D \hookrightarrow \oO_{D_\p}^\ast \to
\rightquom{\oO_{D_\p}^\ast}{1 + \sQ}{3pt}{\big} \cong \FF_{q^2}^\ast
\]
See e.g.~\cite[\S6.4]{MaclachlanReidBook} for details.
Set $\Gamma(\p) \lhd \Gamma$ to be the projectivization of the kernel of
$\oO_D^\ast \to \FF_{q^2}^\ast$, and let $X(\p) = \Gamma(\p) \backslash \hthree$.  The
precise degree of $X(\p) \to X$ depends on $K$ and $D$, as for instance
$\oO^*_D \to \FF_{q^2}$ is rarely onto.

More generally, suppose $A = \p_1 \p_2 \cdots \p_n$ is a square-free
ideal of $\oO_K$.  Then $X(A)$ is defined as the common cover of the
$X(\p_i)$.  If $D$ does not ramify at any $\p_i$, strong
approximation implies that the degree of $X(A) \to X$ is the product
of the degrees of $X(\p_i) \to X$.

\subsection{Standard Hecke operators}\label{subsec-std-hecke-ops}

The Hecke operators for these covers are defined as follows.  Let
$\ram(D)$ denote the set of primes where $D$ is ramified.  Suppose
$A$ divides a square-free ideal $B$. Then we have the natural
covering map
\[
\phi_1=\phi_{1,B/A} \maps X(B) \to X(A).
\]
Moreover, if $B/A$ is prime to ram$(D)$, then for each ideal $J$
dividing $B/A$, there is a certain covering map:
\[
\phi_J=\phi_{J,B/A} \maps X(B) \to X(A),
\]
which is just $\phi_1$ when $J$ is the unit ideal. These maps are
defined below. In particular, for a prime $\p\notin\ram(D)$ not
dividing $A$, we have two covering maps:
\[
\phi_1, \phi_\p \maps X(A\p) \to X(A)
\]
Then the Hecke operator for $\p$ on $H^*(X(A))$ is defined as $T_\p =
T_{\phi_1, \phi_\p}$ in the notation of Section~\ref{sec_hecke_ops}.
These have the following key properties:
\begin{proposition} \label{prop-hecke-props}
  Let $A$ be a square-free ideal of $\oO_K$.
\begin{enumerate}

\item \label{hecke_comm}
  The standard Hecke algebra for $X(A)$ is generated by the $T_\p$ for
  $p \nmid A$ and $p \notin \ram(D)$.  For prime ideals $\p$ and $\q$, the
  Hecke operators $T_\p$ and $T_\q$ commute.

\item \label{comm_degen}
  Consider a cover $\phi_J \maps X(B) \to X(A)$, with $B$ square-free,
  and the associated
  \emph{degeneracy map}
  \[
  \phi_J^* \maps H^i(X(A)) \to H^i(X(B)).
  \]
  If
  $\p$ is a prime ideal not dividing $B$, then
  \[
  \phi_J^* \circ T_\p^A  = T_\p^B \circ \phi_J^*
  \]
  where the superscripts on the Hecke operators indicate which
  orifolds' cohomology is being acted upon.
\end{enumerate}
\end{proposition}

We will now give a definition of the maps $\phi_J$, and indicate how
to deduce the properties above. For this it is convenient to work in
the ad\`elic setup.

Let $\A_K=K_\infty\times\A_{K,f}$ be the addle ring of $K$, with
$K_\infty$ denoting the product of the archimedean completions of
$K$ and $\A_{K,f}$ being the ring of finite ad\`eles. Put $G=D^\ast$,
which is a reductive algebraic group over $K$, with center $Z$. Then
$G(K)$ is a discrete subgroup of the locally compact group $G(\A_K)
= G_\infty\times G(\A_{K,f})$. Since $K$ has a unique complex place
and $D$ is a quaternion division algebra ramified at all the real
places, we have $G_\infty \simeq \GL{2}{\C}\times {\mathbb
H}^{[K:\Q]-2}$, where ${\mathbb H}$ is Hamilton's quaternion algebra
over $\R$; consequently, $Z_\infty\backslash G_\infty$ is
$\PGL{2}{\C}\times \mathrm{SU}(2)^{[K:\Q]-2}$. Moreover, since $G$ is
anisotropic, Godement's compactness criterion implies that the
ad\`elic quotient
\[
X_\A: = \leftquom{G(K)Z(\A_K)}{\hthree \times G\left(\A_{K,f}\right)}{4pt}{\big}
\]
is a compact space, equipped with a right action by $G(\A_{K,f})$.
We may view $X_\A$ as the projective limit of closed hyperbolic
\3-orbifolds
\[
X_U: = X_\A/U,
\]
as $U$ varies over a cofinal system of compact open subgroups of
$G(\A_{K,f})$. These orbifolds are in fact manifolds for deep enough
$U$. The reduced norm on $D^\ast$ induces a map $\Nrd: G(\A_{K,f})\to
\A_{K,f}^\ast$, and if $U$ is any compact open subgroup of
$G(\A_{K,f})$, we can write $G(\A_K)$ as a finite union
\[
G(\A_K)  =  \coprod_{j=1}^{h(U)} G(K)G_\infty x_jU,
\]
for any elements $x_1=1, x_2, \dots, x_{h(U)}$ in $G(\A_{K,f})$ such
that $\{\Nrd(x_j)\}$ is a complete set of representatives for
$K^\ast \Nrd(U)\backslash \A_{K,f}^\ast$.  We may choose these $x_j$ to
have components $1$ at the primes in $\ram(D)$, and such that for
finite $v\notin\ram(D)$,
\[
x_{j,v} \, = \, \begin{pmatrix}a_{j,v} & 0\\0 & 1\end{pmatrix},
\]
for some $a_{j,v}$ in $K_v^\ast$, of course with $a_{j,v}\in
\oO_{K,v}^\ast$ for almost all $v$. When $U$ is the maximal compact
subgroup $G(\oO_{K,f})=\prod_P G(\oO_{K,P})$, the number of components
$h(U)$ is just the class number $h$ of $K$.

We will use the following cofinal system of compact open subgroups
$U_A$ for ideals $A=\prod_P P^{m(P)}$, given by:
\[
U_A =  \left(\prod_{P\in {\ram}(D)} U_P(P^{m_P})\right)\times \left(\prod_{P\notin {\ram}(D)} U_{0,P}(P^{m_P})\right),
\]
where for $P\in \ram(D)$, the subgroup $U_P(P^{m_P})$ is the principal
congruence subgroup of $G(K_P)$ of level $P^{m_P}$, while for
$P\notin \ram(D)$,
\[
U_{0,P}(P^{m_P})  =  \setdef{g\in G(\oO_{K,P})}{g \equiv
\begin{pmatrix}\ast & \ast\\0 & \ast\end{pmatrix} \, \mod \, P^{m_P}}.
\]
Note that $\Nrd(U_{0,P}(P^{m_P}))$ is $\oO_{K,P}^\ast$, implying that
\[
(A,{\ram}(D))=1  \implies  h(A)\coleq h(U_A) = h.
\]
For the main example of this paper $K=\Q[\sqrt{-3}]$ for which $h=1$.

The space $X_{U_A}$ has $h(A)$ connected components. More precisely,
if we put
\[
\Gamma_j(A)  =  G(K) \cap x_jU_Ax_j^{-1},
\]
we get
\[
X_{U_A}  =  \coprod_{j=1}^{h(A)}  \leftquom{\Gamma_j(A)}{\hthree}{2pt}{\big}.
\]
Here $\Gamma_1(A)\backslash \hthree$ is the orbifold $X(A)$ that we defined in
Section~\ref{cong_covers}.

\subsection{Degeneracy maps}

We now turn to the degeneracy maps discussed in
Proposition~\ref{prop-hecke-props}(\ref{comm_degen}).  If $P$ is a
prime ideal not in $\ram(D)$, choose a uniformizer $\varpi$ at $P$ and
define $g(P) \in G(\A_K)$ to be $\mysmallmatrix{1}{0}{0}{\varpi}$ at $P$
and $1$ at all other places.  We then define a map via right
multiplication:
\[
X_\A \to X_\A \quad \mbox{given by $\xi\mapsto \xi g(P)$}.
\]
Note that $g(P)$ has only one non-trivial component, namely at $P$,
and that component is diagonal. At the finite
levels we get the induced maps
\[
\phi'_P \maps  X_{U_A\cap g(P)U_Ag(P)^{-1}}  \to  X_{U_A}.
\]
It is immediate that
\[
P\nmid A \implies  U_A \cap g(P)U_Ag(P)^{-1} = U_{AP}.
\]
We let $\phi'_P$ denote again the restriction to the connected
component $X(AP)$. It is easy to extend this to square-free ideals
$J$ away from $A$ and $\ram(D)$, and if $J$ divides $B/A$ for
another square-free ideal $B$, then we may compose the natural map
$\phi_1 \maps X(B)\to X(AJ)$ with $\phi_J \maps X(AJ)\to X(A)$ to
obtain the map $\phi_J$ discussed above.

For each prime $P\nmid A$, $P\notin \ram(D)$, the Hecke correspondence
$T_P: X(A) \to X(A)$, defined above by pulling back along $\phi_1$ and
then pushing down by $\phi_P$, also acts on the cohomology of $X(A)$ in
the obvious way.  This Hecke operator does not depend on our choice of
uniformizer, because if we replace $\varpi$ by another uniformizer $\varpi'$ at
$P$, then $u\coleq\varpi'\varpi^{-1}$ is a unit at $P$, and right
multiplication by $\mysmallmatrix{1}{0}{0}{u}$ acts
trivially on $X(A)$. So the action of $T_P$ on $H^\ast(X(A))$ is
invariantly defined for any such $P$, and is denoted as $T_P^A$ if we
want to keep track of the level $A$.

Note that right multiplication by $g(P)$ on the ad\`elic space $X_\A$
induces a family of maps
\[
g(P)^A: X_{U_P} \to  X_U,
\]
as $U=U_A$ runs over congruence subgroups of level $A$ prime to
$P$.  Also, the right multiplication by $g(P)$ and $g(Q)$ on $X_\A$
evidently commute for two distinct primes $P, Q$, since the matrices
$g(P)$ and $g(Q)$ themselves commute. One gets the identities, for any
compact open subgroup $U$ of level $A$ prime to $PQ$,
\begin{align*}
U_{PQ}  &=  U\cap\left(g(P)g(Q)U g(Q)^{-1}g(P)^{-1}\right) \\
 &= U\cap\left(g(Q)g(P)U g(P)^{-1}g(Q)^{-1}\right) =  U_{QP}.
\end{align*}
Moreover
\[
U_P\cap \left(g(Q)U_P g(Q)^{-1}\right)  =  U_{PQ} \quad \mbox{and} \quad  U_Q\cap \left(g(P)U_Q g(P)^{-1}\right)  = U_{PQ}.
\]
and one also gets a commutative diagram
\[
\xymatrix @ur {
X(AP) \ar[d]_{g(P)^{AQ}} & X({APQ})  \ar[l]_{g(Q)^{AP}} \ar[d]^{g(P)^{AQ}} \\
X(A) & X(AQ)  \ar[l]^{g(Q)^A}
}
\]
From this, and the definition of $T_P^A$, we get assertion (1) of
Proposition 3.9. For assertion (2) we note that $T_P$ descends
compatibly to levels $A$ and $B$ with $A \mid B$, as long as $P
\nmid B$. Compatibility with $\phi_J^\ast$ follows.

\subsection{Newforms and oldforms}
\label{sec_new_and_old}

The cohomology of $X(A)$ decomposes naturally into oldforms and
newforms.  The former, denoted $H^i(X(A))\oldforms$ is by definition
the subspace of $H^i(X(A))$ generated by the images of
$\varphi_{J,B/A}^*$, as $B$ varies over all the ideals properly
containing $A$ and $J$ varies over divisors of $B/A$ which are prime
to ram$(D)$. Note that an oldform is therefore fixed by a conjugate
of a congruence subgroup of smaller level.

As we saw above, the ad\`elic quotient $X(U_A)$ is not connected, and
we have
\[
H^3(X(U_A)) =  H^3(X_\A)^{U_A}  \simeq  \Q^{h(A)},
\]
and this holds even when $X(U_A)$ is an orbifold. As $X(A)$ is the
connected component, it follows that
\[
H^3(X(A))  \simeq  \Q.
\]

Let us focus on $H^1$. There is a complementary \emph{new subspace}
$H^1(X(A))\newforms$, which consists of cohomology classes which are
genuine to the level at hand, and it can be defined as the annihilator
of $H^2(X(A))\oldforms$ under the cup product pairing
\[
\cup \maps  H^1(X(A)) \times H^2(X(A))  \to  H^3(X(A))
\simeq \Q.
\]
Since the elements of the Hecke algebra act as correspondences on
$X(A)$, this pairing is also functorial for the Hecke action.

\subsection{The main geometric theorem}
\label{sec_cong_cover_thm}

Now we turn to
\begin{theorem}\label{cong_cover_thm}
  Let $M = X(A)$ be an arithmetic hyperbolic \3-manifold defined as
  above from a quaternion algebra $D/K$ and an ideal $A$ of $\oO_K$.
  Let $\p_1, \p_2, \ldots,\p_n$ be prime ideals of $\oO_K$ coprime to $A$ at
  which $D$ does not ramify.  Consider the corresponding congruence
  cover $M(\p_1 \p_2 \cdots \p_n) \to M$, which has degree $\prod
  \left(1 + N_{K/\Q}(\p_i) \right)$.  Suppose $\omega \in H^1(M)$ comes
  from a fibration over the circle, and that $T_{\p_i} (\omega) = 0$ for
  each $\p_i$.  Then $M(\p_1 \p_2 \cdots \p_n)$ has at least $2^n$ pairs
  of genuinely distinct fibered faces.
\end{theorem}

Theorem~\ref{intro_thm} follows from Theorem~\ref{cong_cover_thm}
once we exhibit in Section~\ref{sec_example} an $M$ and an $\omega$,
together with an \emph{infinite} set of primes $\{\p_i\}$ such that
$T_{\p_i}(\omega) = 0$, by considering all the covers
$M(\p_1\p_2\cdots \p_m)$.

\begin{proof}
  We claim inductively that $M_m = M(\p_1 \p_2 \cdots \p_m)$ has $2^m$
  fibered classes $\omega_i$ lying in genuinely distinct fibered faces
  which are killed by $T_{\p_k}$ for all $k > m$.  To simplify
  notation, set $\p = \p_{m + 1}$.  By the discussion above, we have
  two covering maps
  \[
  \phi_1, \phi_\p \maps M_{m + 1} \to M_m
  \]
  and let $\tilde \omega_j$ be the $2^{m + 1}$ classes which are the
  pull-backs of the $\omega_i$ by $\phi_1$ and $\phi_\p$.  As $T_\p$ kills all
  the $\omega_i$, Lemma~\ref{full_geom_thm} implies that the $\tilde \omega_j$
  lie in genuinely distinct fibered faces.  The commutatively property
  of Hecke operators and degeneracy maps implies that $T_{p_k} (\tilde
  \omega_j) = 0$ for all $j$ and any $k > m + 1$, completing the
  induction.
\end{proof}

\section{The base manifold $M$}
\label{sec_example}

Here is the arithmetic description of our example in the framework of
Section~\ref{cong_covers}.  Consider the imaginary quadratic field
$K=\Q[\sqrt{-3}]$, with ring of integers $\oO_K$.  The non-trivial
automorphism of $K$ can be identified with complex conjugation via
embedding $K$ into $\C$, and will be denoted as such.  The rational
prime $7$, or rather the ideal $7 \oO_K$, splits in $K$ as $\q \qbar$.
Let $D$ be the unique quaternion division algebra over $K$ which
ramifies exactly at $\q$ and $\qbar$.  Let $\oO_D$ be a maximal order
of $D$; this order is unique up to conjugation as $K$ is quadratic and
has restricted class number $h_\infty = 1$.  The corresponding hyperbolic
\3-orbifold is $X = \oO^\times_D \backslash \hthree$.  Our example manifold $M$ is
the (principal) congruence cover of $X$ of level $\q \qbar$.  To state
its key properties, we need one more piece of notation.  Let $\sP$ be
the set of rational primes $p\neq 7$ which are inert in $\Q[\sqrt{-7}]$
but are split in $K$, and consider the set of prime ideals of $\oO_K$
given by $\sP_K= \setdef{P}{N_{K/\Q}(P)=p\in \sP}$.  We will show:
\begin{theorem}\label{thm_example_properties}
  Let $M$ be the arithmetic hyperbolic \3-manifold described above.
  \begin{enumerate}
  \item\label{ex_number_part} The new subspace $V = H^1(M; \Q)\newforms$ is
    2-dimensional and isotypic under the Hecke action.  Moreover, for
    each prime $\p$ in the set $\sP_K$ described above, $V$ is in the
    kernel of the Hecke operator $T_\p$.
  \item\label{ex_top_part}
    There is an $\omega \in V$ coming from a fibration of $M$ over the circle.
  \end{enumerate}
\end{theorem}

The proofs of parts (\ref{ex_number_part}) and (\ref{ex_top_part}) are
essentially independent, and we tackle them separately in the next two
sections.  The main result of this paper, Theorem~\ref{rapid_growth},
follows from considering the congruence covers $M(\p_1 \p_2 \cdots \p_n)
\to M$ where $\p_j \in \sP_K$.  Theorem~\ref{thm_example_properties}
implies that $M(\p_1 \p_2 \cdots \p_n)$ satisfies the hypotheses of
Theorem~\ref{cong_cover_thm}, and hence has at least $2^{2n}$ fibered
faces.  Using that $\sP$ consists of $1/4$ of all rational primes, we
then calculate a lower bound on the number of fibered faces in terms
of the degree of the cover.  The details of this proof of
Theorem~\ref{rapid_growth} are given in
Section~\ref{sec_proof_of_main_thm}.

The manifold $M$ is the only one we could find with a CM form coming
from a fibration over the circle.  Because CM forms are fairly rare,
there are very few potential examples where one can computationally
examine the topology in order to check fibering, even with the improved
methods we introduce here.  We were very fortunate that $M$ does
indeed have the desired properties, since there are probably only one or
two more potential examples that are within reach.  The next example to
look at would have been starting with the elliptic curve $y^2+y = x^3$
over $\Q$ which has conductor 27 and is CM by $\Q(\sqrt{-3})$ and then
base-changing to $\Q(\sqrt{-2})$.

\section{The automorphic structure of the cohomology of $M$}
\label{sec-automorphic}

In this section, we give the arithmetic construction of the needed
cohomology class on $M$; the reader is urged to refer back to
Section~\ref{sketch_of_construction} for an outline before proceeding
(the basic construction here has was introduced in
\cite{LabesseSchwermer1986}).  We show in particular why the underlying
automorphic form $\pi_K$ on GL$(2)/K$, $K=\Q[\sqrt{-3}]$, transfers to
the quaternion division algebra $D/K$ ramified only at the primes
above $7$, and moreover determine exactly its level, which is crucial.
The final result of the section gives a specific (closed) arithmetic
hyperbolic $3$-manifold $M$ (of level $7$) such that the new subspace
of $H^1(M,\R)$ contains a plane defined by a CM form $\pi_K^D$ which is
killed by half of the Hecke operators $T_P$.  Later in Section 6, we
will show that this new subspace is at most $2$-dimensional and
contains a vector coming from an $S^1$-fibration.  Taken together,
this will prove Theorem~\ref{thm_example_properties}.

The proofs are given in detail, especially when there are no precisely
quotable references, even though the experts may well be aware of the
various steps.  While there are an abundance of CM forms which define
similar cohomology classes on corresponding manifolds, so far we have
managed to find only this one particular form at the confluence of the
required arithmetic \emph{and} topological properties. Thus we wish to
take extra care in checking the details.

\subsection{A cusp form form on $\GL{2}{\Q}$}

Let $E$ be the elliptic curve over $\Q$ defined by
\[
y^2+xy=x^3-x^2-2x-1.
\]
It has conductor $49$, admitting complex multiplication by an order
in the imaginary quadratic field $L\coleq\Q[\sqrt{-7}]$. For any
rational prime $p$, we set, as usual,
\[
a_p(E) \coleq  p+1-N_p(E),
\]
where $N_p(E)$ is the number of points of $E$ modulo
$p$.  The
$L$-function of $E$ is given by
\[
L(s, E): = \prod_p \left(1-a_p(E)p^{-s}+p^{1-2s}\right)^{-1}
\]
which, by a theorem of Deuring (see e.g.~\cite{Gross}), equals
\[
L(s, \Psi)\coleq \, \prod_{P} \left(1-\Psi({P}){\Norm}_{L/\Q}({P})^{-s}\right)^{-1} =
\prod_p \left(\prod\limits_{{P}\vert p}
\left(1-\Psi({P}){\Norm}_{L/\Q}(P)^{-s}\right)\right)^{-1}
\]
where $\Psi$ is a Hecke character of $L$ of weight $1$, and ${\Norm}_{L/\Q}$
denotes the norm from $L$ to $\Q$. Equating, for each $p$, the
corresponding Euler $p$-factors of $L(s,E)$ and $L(s,\Psi)$, we see that
\[
p\oO_L \quad \mbox{prime} \, \, \implies \, a_p(E) \, = \,  0.
\]

One knows by Hecke, or by using the converse theorem, that $\Psi$
defines a holomorphic cusp form of weight $2$ for the congruence
subgroup $\Gamma_0(49)$ of $\SL{2}{\Z}$, given by $f(z) = \sum_{n\geq 1}
a_nq^n$ with $q=e^{2\pi iz}$ which is new for that level.  Here, we
normalize $f$ so that $a_1(f)=1$ and $a_p(f)$ is the eigenvalue of the
$p^{\mathrm{th}}$ Hecke operator $T_p$.  One obtains, for every prime $p \neq 7$,
that $a_p(f) = a_p(E)$. Thus $a_p(f)$ is zero whenever $p$ is
inert in $L$.

\subsection{Ramification at $7$}

The cusp form $f$ is ramified at $7$, but to make this precise, and to
see what happens when we base change to $K=\Q[\sqrt{-3}]$, it is
necessary for us to consider the situation ad\`elically.

For any number field $k$, let $\A_k$ denote the topological ring of
ad\`eles of $k$, which is a restricted direct product $\prod_v^\prime \Q_v$, as
$v$ runs over the places of $k$ and $k_v$ denotes the completion of
$k$ at $v$.  When $k$ has only one archimedean place up to complex
conjugation, e.g.~$k = \Q$ or $k$ imaginary quadratic, this place is
denoted by $\infty$.  Given any algebraic group over $k$, it makes
sense to consider the locally compact topological group $G(\A_k) =
\prod_v^\prime G(k_v)$.

One knows that the cusp form $f$ generates a unitary, cuspidal
automorphic representation $\pi=\otimes_v^\prime \pi_v$ of $\GL{2}{\A_\Q}$ with
trivial central character, with $\pi_\infty$ being the lowest discrete
series representation of $\GL{2}{\R}$ (\cite{Gelbart}), such that
$L(s-1/2,\pi)=L(s,f)$. By construction, $\pi_p$ is unramified at every
prime $p\neq 7$, which translates to the elliptic curve $E$ having good
reduction at $p$. More precisely, the \emph{conductor} $c(\pi)$ of $\pi$
satisfies
\[
c(\pi) =  c(\pi_7)  =  49.
\]
We need to identify the local representation at $7$, and in fact to
make sure that it is not a principal series representation.

\begin{lemma}\label{ram lem}  Let $f, \pi$ be as above. Then
$\pi_7$ is a supercuspidal representation of $\GL{2}{\Q_7}$ of
conductor $49$.
\end{lemma}

After this paper was written, we learnt from B.~Gross of an earlier,
different proof of this Lemma in \cite{GL}. Moreover, the referee
has remarked that such a result should by now be well known. We will
nevertheless give a proof, as it is explicit, not long, and uses an
explicit criterion of Rohrlich involving the discriminant.

\begin{proof}
  Let us first note a few basic facts. Since $\pi$ is defined by a
  Hecke character of a quadratic extension, or equivalently since $E$
  has complex multiplication, the local representation $\pi_p$ at any
  prime $p$ is either in the principal series or is supercuspidal. In
  either case, we claim that there is a finite extension $F/\Q_p$ over
  which (the base change of) $\pi_p$ becomes an unramified principal
  series representation. Indeed, by the local Langlands correspondence
  for $\GLtwo$ \cite{Kutz}, at any $p$ the local representation $\pi_p$
  is attached to a semisimple $2$-dimensional $\C$-representation
  $\sigma_p$ of the local Weil group $W_{\Q_p}$, whose determinant
  corresponds to the central character of $\pi_p$ (by local class field
  theory), which is trivial in our case. Recall that $W_{\Q_p}$ is the
  dense subgroup of the absolute Galois group $G_p=
  \Gal(\Qbar_p/\Q_p)$ generated by the inertia subgroup and the
  integral powers of a lift to $G_p$ of the Frobenius automorphism
  $x\to x^p$ of $\overline {\FF}_p$. One knows that $\pi_p$ is in the
  principal series if and only if $\sigma_p$ is reducible, necessarily of
  the form $\nu_p\oplus\nu_p^{-1}$. We may write $\nu_p$ as an unramified
  character times a finite order character $\lambda_p$, since the inertia
  subgroup $I_p$ acts by a finite quotient, and so $\sigma_p$ becomes
  unramified over the finite extension $\Q_p(\lambda_p)$, which is the
  cyclic extension of $\Q_p$ cut out by $\lambda_p$ by class field theory.
  By looking at the $\ell$-adic representation of $G_p$ attached to
  $\sigma_p$, and using the N\'eron-Ogg-Shafarevich criterion \cite{Silver},
  we see that $E$ also acquires good reduction over $\Q_p(\lambda_p)$.
  Suppose now $\sigma_p$ is irreducible, in which case $\pi_p$ is
  supercuspidal. Again, $\sigma_p$ is acted on by $I_p$ through a finite
  quotient $H_p$, because $I_p$ is profinite and $\sigma_p$ takes values
  in $\GL{2}{\C}$.  Consequently, there exists a finite Galois
  extension $F/\Q_p$ with Galois group $H_p$ such that over $F$, the
  representation $\sigma_p$ becomes unramified and $E$ attains good
  reduction. Hence the claim.

  We argue further that this line of reasoning shows that $\pi_p$ is in
  the principal series if and only if $E$ acquires good reduction over
  an \emph{abelian} extension of $\Q_p$. Indeed, this is clear if
  $\pi_p$ is in the principal series, so we may assume that $\pi_p$ is
  not of this kind.  Then $\sigma_p$ is irreducible and becomes unramified
  only over a more complicated extension $F$ due to irreducibility;
  $F$ is dihedral for $p \geq 3$, but can be wilder for $p=2$.

  For $p\geq 5$, a theorem of Rohrlich \cite{Rohr} asserts that $E$
  acquires good reduction over an abelian extension of $\Q_p$ if and
  only if the following field is abelian over $\Q_p$:
\[
F\coleq \Q_p(\Delta^{1/e}),
\]
where $\Delta$ is the discriminant of $E$ and
\[
e  =  \frac{12}{(v_p(\Delta),12)},
\]
where the denominator on the right signifies the gcd of
$v_p(\Delta)$ and $12$.

In the present case, $p=7$ and $\Delta=-7^3$, so that $e = 12/3 = 4$.
Consequently, $F$ is generated over $\Q_7$ by a fourth root of $343$.
The only way $F$ can be abelian over $\Q_7$ is for $\Q_7$ to contain
the fourth roots of unity, i.e., for $-1$ to be a square in $\Q_7$,
and hence in $\FF_7$. But this does not happen because the Legendre
symbol $\genfrac{(}{)}{}{}{-1}{7}$ is $(-1)^{(7-1)/2}=-1$. Hence $\pi_7$
is supercuspidal.

Finally, the conductor $c(\pi)$ of $\pi$ factors as
\[
c(\pi)\, = \, \prod_p c(\pi_p),
\]
and since $\pi$ is associated to $E$, the conductor $c(\pi)$ coincides
with that of $E$, which is $49$. Moreover, as $\pi_p$ is unramified at
every prime outside $7$, we have $c(\pi_p)=1$ for every $p\neq 7$. It follows
that $c(\pi_7)=49$.
\end{proof}

\subsection{Local transfer to the quaternion algebra over $\Q_7$}

Denote by $B$ the unique quaternion division algebra over $\Q_7$.
Then $B^\ast$ is an inner form of $\GLtwo$ over $\Q_7$. By the local
Jacquet-Langlands correspondence \cite[\S16]{JL}, there is a
finite-dimensional, irreducible representation $\pi_7^B$ of $B^\ast$
functorially associated to $\pi_7$. Let $\oO_B$ denote the maximal
order of $B$, and for any $m\geq 0$, let $\Gamma_B(7^m)$ denote the
principal congruence subgroup of level $7^m$, with $\Gamma_B(1)$
denoting $\oO_B^\ast$.

\begin{lemma} \label{lem-local-rep}
  The dimension of $\pi_7^B$ is $2$, and its conductor is $49$. This
  representation is non-trivial when restricted to $\oO_B^\ast$, and the
  kernel contains $\Gamma_B(7)$.
\end{lemma}

Some people would be tempted to say that $\pi_7^B$ has conductor
$7$, but it would not be in agreement with the convention used in
the theory of automorphic forms, where the trivial representation of
$B^\ast$ is said to have (normalized) conductor $7$ (see \cite[\S3]{Tu}).

\begin{proof}
Two of the basic properties of the correspondence $\pi_7 \to \pi^B_7$
are the following:
\begin{enumerate}
\item \label{item-formal-dim}
The dimension of $\pi^B_7$ is the formal dimension of $\pi_7$.
\item \label{item-conductor} $c(\pi_7^B) = c(\pi_7)$.
\end{enumerate}
To elaborate, under the Jacquet-Langlands correspondence (see
\cite{GelbartJacquet}, Theorem 8.1, for example), the character of
$\pi_7^B$ is the negative of the (generalized) character
$\mathit{ch}(\pi_7)$ of (the discrete series representation) $\pi_7$
on the respective sets of regular elliptic elements (which are in
bijection); $\mathit{ch}(\pi_7)$ is \emph{\`a
  priori} a distribution, but is represented by a function on the
regular elliptic set which is locally constant. For a proof of the
identity (\ref{item-formal-dim}), which follows from this character
relation, see for example \cite[Prop.~5.9]{Ro}. The point is that
$\mathit{ch}(\pi_7)$ on regular elliptic elements close to $1$
equals, up to sign, its formal dimension, and similarly for the
character of of $\pi_7^B$.  For a discussion of the identity
(\ref{item-conductor}), see e.g.~\cite[page~21]{P}. One can also
refer to \cite[\S56]{BuH}, once one observes that (in the notation
of that book) $c(\pi_7^B)=7^{\ell(\pi_7^B)+1}$ and
$c(\pi_7)=7^{2\ell(\pi_7)+1}$.

Thanks to (\ref{item-formal-dim}), it suffices to check that the
formal dimension of $\pi_7$ is $2$. Let us first make the following

\begin{claim}\label{claim-rep-is-induced}
  The $2$-dimensional representation $\sigma_7$ of $W_{\Q_7}$ associated
  to $\pi_7$ is induced by a character $\chi$ of $W_F$, where $F$ is the
  unique unramified quadratic extension of $\Q_7$.
\end{claim}

To begin, as $\pi$ is associated to a Hecke character $\Psi$ of
$\Q[\sqrt{-7}]$, the representation $\sigma_7$ is induced by a character
$\psi$, of $W_k$, where $k$ is the ramified quadratic extension of
$\Q_7$ obtained by completing $L$ at the prime $(\sqrt{-7})$. (The
character $\psi$ is the local component at $\sqrt{-7}$ of the global
unitary id\`ele class character of $L$ defined by $\Psi$.) Then we have
\[
c(\sigma_7) = \Norm_{k/\Q_7}(c(\psi))\disc_{k/\Q_7},
\]
where $\Norm$ denotes the norm and $\disc$ the discriminant.
Since $-7\equiv 1\pmod 4$,
${\disc}_{k/\Q_7}$ is $-7$, which forces $c(\psi)$ to be $P \coleq
(\sqrt{-7})\oO_k$, because $c(\sigma_7)=49$ and
${\Norm}_{k/\Q_7}(\sqrt{-7})=7$. Moreover, we have
\[
{\det}(\sigma_7)  = \mathit{Ver}_{k/\Q_7}(\psi)\eta_k,
\]
where $\mathit{Ver}$ denotes the \emph{transfer map}, and $\eta_k$ is
the quadratic character $W_{\Q_7}$ attached to $k/\Q_7$. As the
determinant of $\sigma_7$ is $1$, the restriction of $\psi$ to
$\Q_7^\ast$ must be the quadratic character $\eta_k$. If $\tau$
denotes the non-trivial automorphism of $k/\Q_7$, we must have
\[
\psi^\tau  =  \psi\nu,
\]
where $\nu$ is the unique quadratic unramified character of $W_k$.
(As usual, $\psi^\tau(x)= \psi(x^\tau)$, for all $x\in W_k$.)  Indeed, with
$\{U_k^{j}\}$ denting the usual filtration of the unit group $U_k$,
$\psi$ and $\psi^\tau$ are trivial on $U_k^1$, but not on $U_k$, since
$c(\psi)=P$. Also, $\psi/\psi^\tau$ is evidently trivial on $\Q_7^\ast$,
implying that $\psi/\psi^\tau$ is even trivial on the full unit group $U_k$,
and thus must be an unramified character $\nu$. Also, since
$(\psi^\tau)^\tau=\psi$, we get $\nu^2=1$. On the other hand, since
$k=\Q_7(\sqrt{-7})$, the automorphism $\tau$ sends the uniformizer
$\sqrt{-7}$ to its negative, and so $\psi^\tau$ differs from $\psi$; in
other words, $\nu$ must be the unique unramified quadratic character of
$W_k$. Next observe that since $k/\Q_7$ is ramified, the unique
unramified quadratic extension of $k$ is the compositum $Fk$, where
$F$ is the unique unramified quadratic extension of $\Q_7$. This
forces $\nu$ to be of the form $\eta_F\circ {\Norm}_{K/k}$, and so we have
\[
\sigma_7\otimes\eta_F   \simeq   I_k^{\Q_7}(\psi\nu)  \simeq  I_k^{\Q_7}(\psi^\tau) \simeq   \sigma_7,
\]
as $\sigma_7$ is induced equally by $\psi$ and $\psi^\tau$. Thus
$\sigma_7$ admits a self-twist under $\eta_F$, or equivalently, that
there is a $W_{\Q_7}$-homomorphism
\[
\eta_F  \hookrightarrow  \sigma_7\otimes \sigma_7  \simeq  1 \oplus {\rm sym^2}(\sigma_7).
\]
Here we have used the triviality of $\det(\sigma_7)$, and
$\mathrm{sym}^2(\sigma_7)$ denotes the (\3-dimensional) symmetric square
representation. Combining this with $\sigma_7\simeq\sigma_7^\lor$, we get a
$W_{\Q_7}$-isomorphism
\[
{\End}(\sigma_7)  \simeq  \sigma_7^{\otimes 2}  \simeq \eta_F \oplus {\underline 1} \oplus \beta,
\]
for a $2$-dimensional representation $\beta$ ($\subset {\rm
sym^2}(\sigma_7)$). Since $\eta_F$ is trivial on $W_F$, we get
\[
{\dim}_{\C}{\Hom}_{W_F}({\underline 1},{\End}(\sigma_7)) \,
\geq \, 2.
\]
So by Schur, the restriction of $\sigma_7$ to $W_F$ is reducible,
and if $\chi$ is a character of $W_F$ appearing in
${\sigma_7}_{|_F}$, Frobenius reciprocity implies that $\sigma_7$ is
induced by $\chi$.  This proves Claim~\ref{claim-rep-is-induced}.

By the local correspondence \cite{JL}, the representation $\pi_7$ is
associated to the character $\chi$ of $F^\ast$ ($\simeq
W_F^{\mathrm{ab}}$), with $F/\Q_7$ unramified. Now we may apply Lemma
5 of \cite{PrasadRamak}, which summarizes the results we want on the
formal degree, etc. In that lemma, $\mathrm{cond}(-)$ means the
\emph{exponent} of $c(-)$, and the quadratic extension $K$ there
corresponds to our (unramified) $F$ here. (It should be noted that the
$K$ in that lemma is taken to be unramified if the representation is
attached to more than one quadratic extension.) Since
$\mathrm{cond}(\chi)=1$ for us, we see that $\pi_7^B$ has dimension
$2$. Also, $c(\pi_7^B)=1$. This implies the last part of the lemma by
the discussion in \cite{Tu}, cf.~the two paragraphs in Section 3
before Theorem 3.6.
\end{proof}

\subsection{Base change to $\Q[\sqrt{-3}]$ and the global transfer}

Put $K=\Q[\sqrt{-3}]$, in which the prime $7$ splits as follows:
\[
7\oO_K  =  Q{\overline Q}, \quad Q=\left(2+(1+i\sqrt{3})/2\right).
\]
Let $\pi_K = \otimes_v^\prime \pi_{K,v}$ denote the base change
\cite{Langlands} of $\pi$ to $K$, which has conductor
$\q^2\qbarsquare$, the reason being that $7$ splits in $K$ as
above, and consequently,
\[
\pi_{K,\q}\simeq \pi_{K,\qbar} \simeq \pi_7.
\]
Also, $\pi_K$ still has trivial central character. More importantly,
since $\pi_K$ has a supercuspidal component, it must be cuspidal.
One can also see this from Deuring's theory since $K$ and $L$ are
disjoint.

Let $D$ be, once again, the quaternion division algebra over $K$,
which is ramified only at $\q$ and ${\qbar}$. Then one has, by the
Jacquet-Langlands correspondence, a corresponding cuspidal
automorphic representation ${\pi_K^D}$ of $D(\A_K)^\ast$ such that
\[
v \not\in \{\q, {\qbar}\}  \implies  \pi_{K,v}^D =
\pi_{K,v},
\]
which makes sense because $D_V$ and $\GL{2}{K_v}$ are the same at any
place $v \neq \q, {\qbar}$. More importantly,
\[
D_Q  \simeq  D_{\overline Q}  \simeq  B,
\]
where $B$ is the quaternion division algebra considered in the
previous section, and by the naturality of the Jacquet-Langlands
correspondence,
\begin{equation}\label{eq-local-reps}
\pi_{K_\q}^D  \simeq  \pi_{K_{\qbar}}^D  \simeq  \pi_7^B.
\end{equation}

\begin{lemma}\label{D-level lem} The conductor of ${\pi_K^D}$ is $49$.
The corresponding principal congruence subgroup $U$ of
$D^\ast(\A_{K,f})$ is of level $7\oO_K=\q{\qbar}$. Moreover, the
space of vectors fixed by $U$ is of dimension $4$.
\end{lemma}

\begin{proof}  If $v$ is any finite place of $K$ outside
$\{\q,\qbar\}$, the local representation $\pi_{K,v}^D\simeq
\pi_{K,v}$ is unramified and hence has trivial conductor. Thus we
have, by applying Lemma~\ref{lem-local-rep},
\[
c\left(\pi_K^D\right)  =  c\left(\pi_{K,\q}^D\right)c\left(\pi_{K,{\qbar}}^D\right)  =
c\left(\pi_7^B\right)^2  =  7^2  =  \left(\q\qbar\right)^2.
\]
By our convention on the conductor, it follows the corresponding
principal congruence subgroup is of level $7 \oO_K$.  The final claim
of the lemma follows from Lemma~\ref{lem-local-rep}, in view of
(\ref{eq-local-reps}) above.
\end{proof}

\subsection{The $\pi^D$-isotypic subspace}

Put
\[
{\mathcal V}(D): =  L^2\left(\leftquom{Z(\A_K)D^\ast}{D(\A_K)^\ast}{3pt}{\big}\right),
\]
where $Z$ is the center ($\simeq K^\ast$) of the algebraic group
$D^\ast$ over $K$. Then ${\mathcal V}(D)$ is a unitary
representation of $D(\A_K)^\ast$ under the right action.
The coset space $Z(\A_K)D^\ast\backslash D(\A_K)^\ast$ is compact,
which follows by Godement's compactness
criterion as $D^\ast$ contains no unipotents $\neq 1$. Consequently,
one has a (Hilbert) direct sum decomposition (as unitary
$D(\A_K)^\ast$-modules)
\[
{\mathcal V}(D)  \simeq  \oplus_{\eta}  m(\eta) {\mathcal V}_\eta,
\]
where $(\eta, {\mathcal V}_\eta)$ runs over irreducible unitary,
admissible representations of the group $D(\A_K)^\ast$, with
multiplicities $m(\eta) \geq 0$.

The Jacquet-Langlands correspondence which embeds ${\mathcal V}(D)$
into the discrete spectrum of the corresponding space for
$\GL{2}{\A_K}$, and since the multiplicities are $1$ on the
$\GLtwo/K$-side \cite{JL}, we get
\[
m(\eta)  \leq 1, \quad  \mbox{for all $\eta$.}
\]

Since $\pi_K^D = \pi_{K,\infty}\otimes \pi^D_{K,f}$ is by
construction a cuspidal automorphic representation of $D(\A_K)^\ast$
of trivial central character, it gives rise to a non-trivial summand
in ${\mathcal V}(D)$. In particular,
\[
m(\pi^D_K)  =  1.
\]

The $\pi_K^D$-isotypic subspace of ${\mathcal V}(D)$, identifiable
with ${\mathcal V}_{\pi_K^D}$, is infinite\hyp dimensional. However,
due to the admissibility of $\pi^D_{K,f}$, given any compact open
subgroup $U$ of $D(\A_{K,f})^\ast$, the space $U$-invariants in
$\pi^D_{K,f}$ is finite dimensional. Applying Lemma~\ref{D-level
lem}, we obtain the following
\begin{lemma}
\[
{\dim}_\C  {\mathcal V}_{\pi^D_{K,f}}^{U_D(7)}  =  4.
\]
\end{lemma}

\subsection{Connected components of the ad\`elic double coset space}

Now put $G=D^\ast/K$ and consider the ad\`elic quotient space
\[
X_\A: = \leftquom{G(K)Z(\A_K)}{\hthree \times
G\left(\A_{K,f}\right)}{4pt}{\big},
\]
introduced in the middle of Section~\ref{subsec-std-hecke-ops}, which
admits a right action by $G(\A_{K,f})$. Also set
\[
A=Q\qbar, \quad {\rm and} \quad X_{U_A} = X_\A/U_A,
\]
where
\[
U_A =  U_Q(Q)\times U_{\qbar}(\qbar)\times \left(\prod_{P\neq
Q,\qbar} G(\oO_{K_P})\right),
\]
where $U_Q(Q)$ is the principal congruence subgroup of $G(K_Q)$ of
level $1$, i.e., the kernel of the map to $G(\F_7)$, and the same for
$U_\qbar(\qbar) \leq G(K_\qbar)$. We will write $U_D(7)$ instead of
$U_A$ to refer to the specific choice of $A$ here.

\begin{lemma}
  $X_{U_A}$ has two connected components, both diffeomorphic to the manifold $M$
  defined in Section~\ref{sec_example}.
\end{lemma}

\begin{proof} The disconnectedness comes because we are working with the
group $D^\ast$ and not its semisimple part
\[
D^1: = \mathrm{Ker}\left(\Nrd: D^\ast \, \to \, K^\ast\right),
\]
where $\Nrd$ denotes the reduced norm.  In our case, we have
\[
X_{U_A} = \raisebox{-4pt}{$D^\ast$} \Big\backslash
\raisebox{4pt}{$\hthree \times \left( \doublequom{Z(\A_{K,
f})}{D(\A_{K,f})^\ast }{U_D(7)}{2pt}{\big} \right)$}
\]
where $D^\ast$ acts diagonally.  Thus the number of connected components
of $M^\A$ is
\begin{equation}\label{conn_comp_form}
\pi_0(M^\A)  = \doublequom{ D^\ast Z(\A_{K, f})}{D(\A_{K,f})^\ast }{U_D(7)}{3pt}{\big}
\end{equation}
Were we working with $D^1$, we would have
\[
\pi_0 = \doublequom{ D^1}{D(\A_{K,f})^1 }{U}{3pt}{\big}
\]
where the right hand quotient is a single element as $D^1$ is dense in
$D(\A_{K,f})^1$ by strong approximation.

Returning to the case at hand, it is straightforward to check that we
can evaluate the right hand side of (\ref{conn_comp_form}) by taking
its image under $\Nrd$.  Since $\Nrd$ surjects $D^\ast$ onto $K^\ast$, as well as
$D(\A_K)^\ast$ onto $\A_K^\ast$, and $Z(\A_K)$ onto ${\A_K^\ast}^2$, we get
\[
\pi_0(M^\A)  =  \doublequom{{\A_{K,f}^\ast}^2K^\ast}{\A_{K,f}^\ast}{\Nrd(U_D(7))}{3pt}{\big}.
\]

Note that since the class number of $K=\Q[\sqrt{-3}]$ is $1$, we
have by strong approximation,
\[
\A_K^\ast  =  K^\ast K_\infty^\ast\prod_P U_{P},
\]
where $P$ runs over all the (finite) primes of $\oO_K$, and
$U_{P}=\oO_{K_P}^\ast$. Moreover,
\[
{\Nrd}(U_D(7))  =  \left(\prod_{P \neq \q,\qbar}
U_P\right) \times U_\q^1\times U_{\qbar}^1,
\]
where $U_\q^1=1+\q\oO_{K_\q}$.

As $7$ splits into $\q, \qbar$ in $K$, we have $\oO_{K_\q} \simeq
\oO_{K_\qbar} \simeq \Z_7$, and moreover, $\Z_7^\ast/{\Z_7^\ast}^2$
identifies with $\FF_7^\ast/{\FF_7^\ast}^2 \simeq \{\pm 1\}$. It
follows that
\[
\pi_0(M^\A) \ =  \mathrm{Coker}\left(\oO_K^\ast \to
\left(\oO_{K_\q}^\ast/{\oO_{K_\q}^\ast}^2\right)\times\left(\oO_{K_\qbar}^\ast/{\oO_{K_\qbar}^\ast}^2\right)
\simeq \left(\FF_7^\ast/{\FF_7^\ast}^2\right)^2\right).
\]
The map $\oO_K^\ast\to \, \FF_7^\ast \times \FF_7^\ast$ is onto the diagonal,
and thus $M^\A$ has two connected components.

To complete the lemma, here is an explicit description of the two
components.  Let $x=(x_v)$ be an element of $\A_K^\ast$ with $x_v=1$ for
all $v\neq Q$, such that the $Q$-component $x_Q$ is an element of
$\oO_{K_Q}^\ast$ ($\simeq \Z_7^\ast$) mapping onto a non-square of $\FF_7^\ast$.
By the preceding, we get an identification
\[
M^\A  =  M_\Gamma  \coprod   M_{\Gamma^\prime},
\]
with connected components
\[
M_\Gamma=\Gamma\backslash\hthree, \quad M_{\Gamma^\prime}=\Gamma^\prime\backslash\hthree,
\]
where
\[
\Gamma=\Gamma(\q\qbar)= D^\ast\cap U_D(7), \quad
\Gamma^\prime = D^\ast \cap \left(xU_D(7)x^{-1}\right).
\]
Now $M_\Gamma$ is precisely the manifold described in
Section~\ref{sec_example}, as claimed, and $M_{\Gamma'}$ is diffeomorphic
to $M_\Gamma$ since $U_D(7)$ is a normal subgroup of $U_D$.
\end{proof}

\subsection{Cohomology on the compact side}

Put $G = D^\ast/Z$. Then $\Gamma=\Gamma(\q\qbar)$ defines a
torsion-free lattice in $G_\infty=\PGL{2}{\C}$. Since $\hthree$ is
contractible, $M_\Gamma$ is an Eilenberg-MacLane space for $\Gamma$,
and thus
\[
H^\ast(M_\Gamma,\Q)  \simeq  H^\ast(\Gamma,\Q).
\]
Moreover, since $\Gamma$ is cocompact, by a refinement of Shapiro's
lemma \cite{Blanc}, we have
\[
H^\ast(\Gamma, \C)  \simeq  H^\ast_\cont \left(G_\infty, L^2\left(\leftquom{\Gamma}{G_\infty}{1pt}{}\right)^\infty\right),
\]
where $H^\ast_\cont$ denotes continuous cohomology, and the
superscript $\infty$ on (the right regular representation)
$L^2(\Gamma\backslash G_\infty)$, which is the unitary analog of
the group algebra of a finite group, denotes taking the smooth
vectors. Consequently, if $L^2(\Gamma\backslash G_\infty)$
decomposes as a unitary $G_\infty$-module as $\oplus m_\infty(\beta)
W_\beta$,
\[
H^\ast(M_\Gamma, \C)  \simeq  H^\ast_{\cont}\left(G_\infty,
L^2\left(\leftquom{\Gamma}{G_\infty}{1pt}{}\right)^\infty\right)
\simeq  \oplus_{(\beta,W_\beta)} m_\infty(\beta)
H^\ast_{\cont}(G_\infty, W_\beta^\infty).
\]
One knows that there is a unique irreducible, unitary representation
$(\beta_0, W_{\beta_0})$ of $G_\infty$ for which
\[
H^i_{\cont}(G_\infty, W_{\beta_0}^\infty) \, \neq \, 0, \quad \mbox{for $i=1,2$,}
\]
and the dimension is $1$. In fact (see \cite{Clozel1987}, for
example), the $2$-dimensional representation of $W_\C=\C^\ast$
attached to this $\beta_0$ (by the local archimedean correspondence)
is given by
\[
\sigma_\infty: z  \to  \left(\frac{z}{\vert z\vert}\right)
\oplus \left(\frac{\overline z}{\vert z\vert}\right).
\]
By construction, the parameter of our $\pi_{K,\infty}$ ($\simeq
\pi^D_{K,\infty}$) is also of this form. Hence we have
\[
\pi^D_{K,\infty} \simeq  W_{\beta_0}^\infty.
\]

If $L^2(G(K)\backslash G(\A_K)/U_D(7)) \simeq\oplus_\Pi \,
m(\Pi)\Pi$ (as unitary $G(\A_K)$-modules), a straightforward ad\`elic
refinement of the above yields the following identification:
\begin{align*}
H^\ast(M^\A,\C) \, &\simeq \, H^\ast_{\cont}\left(G_\infty,
L^2\left(\doublequom{G(K)}{G(\A_K)}{U_D(7)}{2pt}{\big} \right)\right) \\
 &\simeq  \oplus_\Pi H^\ast_{\cont}(G_\infty,{\mathcal V}_{\Pi_\infty})\otimes\Pi_f^{U_D(7)},
\end{align*}
where ${\mathcal V}_\Pi = {\mathcal V}_{\Pi_\infty}\otimes{\mathcal
V}_{\Pi_f}$ denotes the admissible
subspace of (the space of) $\Pi=\Pi_\infty\otimes\Pi_f$.

Next recall that $m_\Pi$ is $1$, and that $H^i_{\cont}(G_\infty,
W_{\beta_0}^\infty)$ is $1$-dimensional for $i=1,2$. So we get an
isotypic decomposition (for $i \in \{1,2\}$),
\[
H^i(M^\A,\C) \, \simeq \, \oplus_{\{\Pi\, : \, \Pi_\infty\simeq
\pi^D_{K,\infty}\}} \, H^i(M^\A,\Pi_f),
\]
where
\[
{\dim}_\C  H^i(M^\A,\Pi_f)  =  {\dim}_\C {\mathcal V}_{\Pi_f}^{U_D(7)}.
\]
In particular, since our $\pi_K^D$ is one such $\Pi$, and since the
space of $U_D(7)$\hyp invariants of $\pi^D_{K,f}$ is $4$-dimensional, we
obtain

\begin{lemma}
For $i=1,2$, we have ${\dim}_\C H^i(M^\A, \pi^D_{K,f}) =  4.$
\end{lemma}

Combining this with the discussion of the connected components of
$M^\A$ in the previous section, since $M^\A$ has two isomorphic
components, we get the following, where $\sP_K$ is the set of primes
defined in Section~\ref{sec_example}:

\begin{theorem}\label{thm_autom_prop}
  For $i=1,2$, we have $\dim_\C  H^i(M_\Gamma;\C)\newforms \geq 2$.
  More precisely, $H^i(M_\Gamma;\C)\newforms$ contains a plane
  (defined by $\pi^D_{K,f}$) on which the Hecke operators $T_P$ act by zero for all $P \in \sP_K$.
\end{theorem}

\section{The Thurston norm and fibered faces of $M$}
\label{sec-topological-props}

In this section, we study the topological properties of $M$ and show
\begin{theorem}\label{thm_top_prop}
  Let $M = X(\q \qbar)$ be the hyperbolic 3-manifold described in
  Section~\ref{sec_example}.  Then $H^1(M; \Q)$ is 3-dimensional and
  \begin{enumerate}
  \item \label{item_thurston_norm}
    The Thurston norm ball $B_T \subset H^1(M; \R)$ is a
    parallelepiped, i.e.~an affine cube.
  \item \label{item_which_fiber}
    Exactly four of the six faces of $B_T$ are are fibered.
  \item \label{item_oldform_loc}
    The subspace of oldforms $H^1(M; \Q)\oldforms$ contains the line passing
    the barycenters of the non-fibered faces.
  \end{enumerate}
\end{theorem}
\noindent Here, the \emph{barycenter} of a face is the average of
its vertices.  We will now show that if Theorem~\ref{thm_top_prop}
holds, then so does Theorem~\ref{thm_example_properties}.

\begin{proof}[Proof of Theorem~\ref{thm_example_properties} modulo
  Theorem~\ref{thm_top_prop}] By Theorem~\ref{thm_autom_prop}, the
  subspace $V = H^1(M; \Q)\newforms$ contains a 2-dimensional subspace
  $V_0$ on which the Hecke operators $T_P$ act by zero for $P \in
  \sP_K$.  By Theorem~\ref{thm_top_prop}, we know$H^1(M; \Q)$ is
  \3-dimensional and $H^1(M; \Q)\oldforms$ is non-zero; this forces
  $V=V_0$. It remains to show that $V$ intersects the interiors of one
  of the four fibered faces.  We can change coordinates on $H^1(M;
  \Q)$ by an element of $\GL{3}{\Q}$ so that $B_T$ is the standard
  cube with vertices $(\pm1, \pm1, \pm1)$, and the non-fibered faces are the
  ones intersecting the $z$-axis.  By part (\ref{item_oldform_loc}),
  the set of oldforms $W$ is then the $z$-axis itself.  Now $\dim V =
  2$, and among all 2-dimensional subspaces, only two miss the
  interiors of the fibered faces, namely the ones containing a
  vertical edge of $B_T$ where two fibered faces meet.  But both of
  these subspaces contain the $z$-axis of oldforms $W$ and so can't be
  $V$. So there is a fibered class in $V$ as desired.
\end{proof}

The rest of this section is devoted to {\it proving}
Theorem~\ref{thm_top_prop}.  We will produce an explicit
triangulation of $M$ and use a series of tricks to compute the
Thurston norm.  We believe these tricks are new; they are remarkably
effective in many examples, and (we hope) are of independent
interest. Our computations used an extensive variety of software:
\cite{SnapPea, Snap, MAGMA214, t3m, SAGE}. Source code, data files, and
complete computational details are available at \cite{paperwebsite}.

\begin{figure}
\centering \leavevmode
\begin{xyoverpic*}{(115,72)}{scale=2.1}{pictures/link}
,(3,44)*+!DR{(5,4)}
,(96,58)*+!DL{(2,0)}
,(75,41)*+++!D{(0,2)}
\end{xyoverpic*}
\caption{A Dehn surgery description of the hyperbolic \3-orbifold $B$.
Our orientation conventions are given by
{\protect\raisebox{-2pt}{\protect\includegraphics[scale=0.35, angle=90]{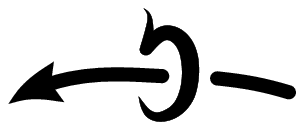}}},
and match those of SnapPea \cite{SnapPea}.}
\label{B-pic}
\end{figure}

\subsection{Finding a topological description}

First, we give a concrete topological description of $M$. Let's
review the setup of Section~\ref{sec_example}.  Let $K =
\Q(\sqrt{-3})$ and let $\q$ and $\qbar$ be the prime ideals sitting
over $7$.  Let $D$ be the quaternion algebra over $K$ ramified
precisely at $\q$ and $\qbar$.  Let $\oO_D$ be a maximal order of
$D$. Let $\phi$ be the embedding of $D^\ast/K^\ast$ into
$\PGL{2}{\C}$, and set $\Gamma = \phi(\oO^\ast_D)$ and $X = \Gamma
\backslash \hthree$.  Also, let $\oO^1_D$ denote the elements of
norm 1, and set $\Gamma' = \phi(\oO^1_D)$ and $X' = \Gamma'
\backslash \hthree$.  We are interested in the congruence covers
$X(\q)$, $X(\qbar)$, and $M = X(\q\qbar)$.

\subsection{The topology of $X'$.}

It is often easier to find $X'$ than $X$ because of the following
characterization of the former.  For a lattice $\Lambda$ in
$\PSL{2}{\C}$, let $\widetilde \Lambda$ denote the preimage lattice
in $\SL{2}{\C}$.  Let $k \Lambda$ be the number field generated by
the traces of elements of $\widetilde \Lambda$, and $A \Lambda$ be
the quaternion algebra generated by the elements of $\widetilde
\Lambda$ (see e.g.~\cite[Ch~3]{MaclachlanReidBook} for details).
Then Corollary~8.3.6 and Theorem~11.1.3 of \cite{MaclachlanReidBook}
imply that $\Lambda \backslash \hthree$ is $X'$ if:
\begin{enumerate}
  \item $k \Lambda = K$ and $A \Lambda = D$.
  \item The trace of every element in $\widetilde \Lambda$ is an
    algebraic integer.
  \item The volume of $\Lambda \backslash \hthree$ is $3^{3/2} 36 \zeta_K(2)/(4 \pi^2) \approx 6.0896496384579219$.
\end{enumerate}
Here, the degree of $K$ is small, as is the volume, so we can expect
these conditions to be verified by the program Snap \cite{SnapPaper,
  Snap} if it is given a topological description of the orbifold
$X'$.

In our case, a brute-force search finds the following description of
$X'$.  Let $B$ be the orbifold shown in Figure~\ref{B-pic}.
From this topological description, SnapPea \cite{SnapPea} derives the
following presentation for $\pi_1(B)$:
\[
\spandef{a,b}{b^2 = aBABaba^4baBAbabABA^4BAb =  aBABaba^3baBABaba^3b = 1},
\]
where $A = a^{-1}$ and $B = b^{-1}$.  Hence $H_1(B; \Z) = \Z/2 \oplus
\Z/8$, and so $B$ has a unique regular cover $C$ with covering group
$\Z/2 \oplus \Z/2$.  SnapPea can explicitly build $C$ as a cover of $B$ given the
action of $\pi_1(B)$ on the cosets of $\pi_1(C)$. Snap checks that
$C$ is $X'$.  This is a little subtle as the reasonable spun ideal
triangulations of $C$ have some flat ideal tetrahedra (though no
negatively oriented ones), and it is necessary to invoke
\cite{PetronioWeeks2000}, in slightly modified form, to justify that
we have constructed the hyperbolic structure.  For details, see
\cite{paperwebsite}.

Though we did not check this rigorously, and will not use it in this
paper, the manifold $X$ is almost certainly the 2-fold cover of $B$
whose homology is $\Z/8 \oplus \Z/8$; this cover is the SnapPea
orbifold $s883(2,-2)(2,-2)$.

\begin{remark}
  Prior to 2007/9/12, the version of the SnapPea kernel \cite{SnapPea}
  that came with SnapPeaPython \cite{SnapPeaPython} contained a bug
  where it sometimes returned the wrong orbifold when the cover it was
  trying to construct was a proper orbifold.  The underlying topology
  was correct but the orbifold loci are sometimes mislabeled by proper
  divisors of their real orders.  This bug also affects the MacOS
  Classic versions of SnapPea, e.g.~SnapPea 2.6 and the derived
  Windows port.  The calculations here were, of course, done with a
  corrected version of the kernel.
\end{remark}

\subsection{Constructing $X(\q)$}

In order to build $X(\q)$, we begin by claiming that it covers $X'$
and that this cover $X(\q) \to X'$ is a cyclic cover of degree 4.
First, we work out the structure of the local quaternion algebra
$D_\q = K_\q \otimes_K D$ by following
\cite[\S6.4]{MaclachlanReidBook}.  We have $K_\q \cong \Q_7$ by
Hensel's lemma since $-3$ has a square-root mod 7.  Let $F$ be the
unique unramified quadratic extension of $K_\q$, concretely $F =
\Q_7(i)$.  Then we have $D_\q \cong F \oplus F j$ where $i$ and $j$
anti-commute and $j^2 = 7$. Then $\oO_{D_\q} \cong \oO_F \oplus
\oO_F j$. Consider the homomorphism
\[
\oO_{D_\q}^\ast  \to \overline{F}^\ast \mtext{defined by $x + y j \mapsto \overline{x}$,
  where $\overline{F}\cong\FF_{49}$ is the residue field of $F$.}
\]
If $\Lambda$ is the kernel of the restricted homomorphism
$\oO^\ast_D \to \overline{F}^*$, then by definition $X(\q) =
\Lambda \backslash \hthree$.

Now the norm map $\Nrd \maps D_\q \to K_\q$ is given by $\Nrd(x + y
j) = \Nrd(x) - 7\Nrd(y)$.  Thus if $\alpha$ is in the kernel of
$\oO^\ast_D \to \overline{F}^\ast$, we have $\Nrd(\alpha) = 1$ in
$\overline{K}_\q$. As the units of $\oO_K$, the sixth roots of
unity, all map to distinct elements in $\overline{K}_\q^\ast$, it
follows that $\alpha$ is in fact in $\oO^1_D$. Thus we have $\Lambda
\subset \oO^1_D$ and so $X'(\q) = X(\q)$; henceforth we focus on the
former description.

The image of $\oO^1_{D_\q} \to \overline{F}^\times$ is just the elements of
$\FF_{49}/ \FF_7$-norm $1$; these elements form a (cyclic) subgroup
$C$ of order 8.  By strong approximation, the image of $\oO^1 \to
\overline{F}^\times$ is the same.  Now $-1 \notin \Lambda$, and so $\Lambda \cong
\pi_1(X'(\q))$.  On the other hand, $-1$ maps into $C$ non-trivially,
and so we conclude that $\pi_1(X'(\q))$ is a normal subgroup of
$\pi_1(X')$ with quotient $C/\{\pm1\} \cong C_4$.

Further, we can characterize $\Lambda$ as exactly those $\alpha \in \oO^1_D$
such that $\tr(\alpha) \equiv 2 \mod \q$.  This is because if $x \in
\overline{F}$ has norm $1$ and $\tr(x) = 2$ then $x = 1$ because $x$
satisfies $x^2 - \tr(x) x + \Nrd(x)$.  Searching through the normal
subgroups of $\pi_1(X')$ with quotient $C_4$, one quickly finds a
unique subgroup for which all tested elements have trace which is $\pm2
\pmod \q$ in the topologically given $\PSL{2}{\C}$-representation; this
subgroup must be $\Lambda$.  In fact, $X(\q)$ and $X(\qbar)$ are
exactly those 4-fold cyclic covers of $X'$ whose homology is $\Z/28 \oplus
\Z/28 \oplus \Z$.
Repeating the same procedure for $\qbar$ builds $\pi_1(X(\qbar))$ as a
subgroup of $\pi_1(X') \leq \pi_1(B)$.  Taking their intersection yields
$\pi_1\big(X(\q\qbar)\big)$.  From the associated permutation action on
$\pi_1(B)\big/\pi_1\big(X(\q\qbar)\big)$, SnapPea can build an explicit
triangulation for $M = X(\q\qbar)$.

\subsection{The Thurston norm of $M = X(\q\qbar)$}\label{subsec_ex_thurston}
For this manifold $M$, we calculate $H_1(M;\Z) = \Z^3 \oplus (\Z/14)^4 \oplus
(\Z/2)^2$.  Our next task is to compute the Thurston norm ball $B_T$
in $H^1(M ; \R) = \R^3$.  In general, there is an algorithm using
normal surfaces for computing the Thurston norm.  Normal surfaces are
those which cut through each simplex in the triangulation like a union
of affine planes, and have been the bedrock algorithmic tool in
\3-dimensional topology since Haken first used them to detect the unknot
\cite{Haken1961}.   For the Thurston norm,  the running time of
this normal surface algorithm is ``only'' simply exponential in the size of the
triangulation \cite{CooperTillmann2007}, but any triangulation of $M$
must have at least 95 tetrahedra as its volume is $\approx 97.434394$, and
the triangulations we constructed had $130$ or more.  At that
complexity, normal surface methods appear to be hopeless.  Instead, we
take a different approach, which involves producing normal surfaces
representing certain classes in $H_2(M)$ quite cheaply.

The other key tool will be the Alexander polynomial $\Delta_M$ of $M$
which lies in the group ring
$\Z\left[H_1(M;\Z)/(\mbox{torsion})\right]$.   Following McMullen
\cite{McMullenNorm},  if $\Delta_M = \sum a_i g_i$ for $a_i \in \Z$ and
$g_i \in H_1(M; \Z)$ then we define the \emph{Alexander norm} on $H^1(M;\R)$  by
\[
\| \omega \|_A  = \sup_{i,j} \omega(g_i - g_j).
\]
The unit ball $B_A$ of this norm is the dual polyhedron to the
Newton polytope of $\Delta_M$ inside $H_1(M;\R)$, which is the
convex hull of the $g_i$.  Since $b_1(M) = \dim H^1(M; \R) > 1$, we
have $\| \omega \|_A \leq \| \omega \|_T$ for all $\omega \in H^1(M;
\R)$, or equivalently $B_T \subset B_A$ \cite{McMullenNorm}. For a
general \3-manifold, these norms do not always coincide, but it is
actually quite common for them to do so, and we will show
\begin{lemma}\label{lemma-norms-equal}
  The Alexander and Thurston norms agree for $M = X(\q \qbar)$.
\end{lemma}

Computing directly from a presentation of $\pi_1(M)$ yields:
\begin{equation}\label{eq_alex_poly}
\Delta_M = \left( 16 x y z  - x y - x z - y  - z  + 16\right)^4
\end{equation}
 Taking
powers of a polynomial just dilates its Newton polytope, so the
Newton polytope of $\Delta_M$ is, up the change of basis, an
octahedron. Dually, this means that $B_A$ is a cube. Since $B_T$ is
a convex subset of $B_A$, to prove Lemma~\ref{lemma-norms-equal} it
suffices to check that the two norms agree on the vertices of the
cube. Moreover, $M$ is very symmetric as $\pi_1(M) \triangleleft
\pi_1(B)$, and any symmetry of $M$ preserves both norms, as does the
map $\iota \maps \omega \to -\omega$ on $H^1(M; \R)$.  One
calculates that the image of
\[
\left\langle \iota, \pi_1(B)  \right\rangle \to \Aut\left(H^1(M; \R)\right)
\]
has order 16.  Now the full symmetry group of the cube has order 48,
and any subgroup of order 16 acts transitively on the vertices, since
the full vertex stabilizer has order 6.  Thus, it suffices to check
$\| \omega \|_A = \| \omega \|_T$ for a single vertex of the cube.

For computing the Thurston norm of a single element of $\omega \in H^1(M;
\Z)$ we used the following method.  Choose a (non-classical)
triangulation $\mathcal T$ of $M$ which has only one vertex, \`a la
Jaco-Rubinstein \cite{JacoRubinstein0eff}.  Then there is a unique
simplicial 1-cocycle representing the given class $\omega$.  As observed
in \cite{Calegari:normal}, any 1-cocycle gives a canonical map $M \to S^1 =
\R/\Z$ which is affine on each simplex of $\mathcal T$.  More
precisely, the integer $\omega$ assigns to an edge specifies how many
times to wrap it around the circle, and this extends over each
simplex because of the cocycle condition.  The inverse image of a
generic point in $S^1$ is a normal surface which is Poincar\'e dual
to $\omega$.  By randomizing the triangulation $\T$ a few times using
Pachner moves, we were able to find a genus 3 surface representing one
of the vertices of $B_A$, which shows the two norms agree for that
class.  This proves Lemma~\ref{lemma-norms-equal}, and thus part
(\ref{item_thurston_norm}) of Theorem~\ref{thm_top_prop}.

\begin{remark}
  In practice, this seems to be a very effective method for computing
  the Thurston norm when the first betti number $b_1$ is $1$, even when
  the triangulation has a couple hundred tetrahedra.  When $b_1 > 1$, things
  are more difficult, as it seems one cannot expect to find a
  triangulation where these special dual surfaces realize the Thurston norm
  for all $\omega$.  While here we evaded this issue by using symmetry, a
  general approach would be to carry along cohomology isomorphisms
  between the new triangulations and the original one while changing
  the triangulation using Pachner moves.  This should enable
  probing the whole Thurston norm ball using this method.
\end{remark}

\subsection{Fibering of $M$}\label{subsec-fibering-M}

Next, we determine which faces of $B_T$ are fibered.  Again, there is
a normal surface based algorithm to decide this \cite{Saul:thesis,
  TollefsonWang1996, JacoTollefson1995}. However, instead of using it,
we build on the method described above.  First, we use the
symmetries of $M$ to simplify the problem.  We know each face of
$B_T = B_A$ is associated to a vertex of the Newton polytope of
$\Delta_M$, and hence has a coefficient of $\Delta_M$ associated to
it.  A generalization of the classical Alexander polynomial
obstruction to fibering says that if a face of $B_T$ is fibered,
then this coefficient must be $\pm1$ (see
e.g.~\cite[Thm.~5.1]{Dunfield:norms}).  In our case, these vertex
coefficients are $\{ 16^4, 16^4, 1, 1, 1, 1\}$; recall from
(\ref{eq_alex_poly}) that $\Delta_M$ is a power of a simpler
polynomial $f$, and so the coefficients are the corresponding powers
of the vertex coefficients of the Newton polytope of $f$. Thus there is
one pair of faces which are definitely not fibered (since their
coefficients are $16^4$), and two pairs of faces which \emph{likely}
fiber. Recall we have a group $S$ of order 16 acting on $H^1(M;\R)$
which preserves $B_T = B_A$, including the coefficient labels on the
faces (up to sign). Thus $S$ preserves the two faces with labels
$16^4$ and has a subgroup of order 8 which leaves each of them
invariant; this subgroup acts there by the full symmetry group of
the square, the dihedral group of order 8.  In particular, $S$ acts
transitively on the remaining four faces of $B_A$.  Thus to show
that all these faces fiber, and so prove part
(\ref{item_which_fiber}) of Theorem~\ref{thm_top_prop}, it suffices
to show:

\begin{lemma}
  The manifold $M = X(\q \qbar)$ fibers over the circle.
\end{lemma}

\begin{figure}
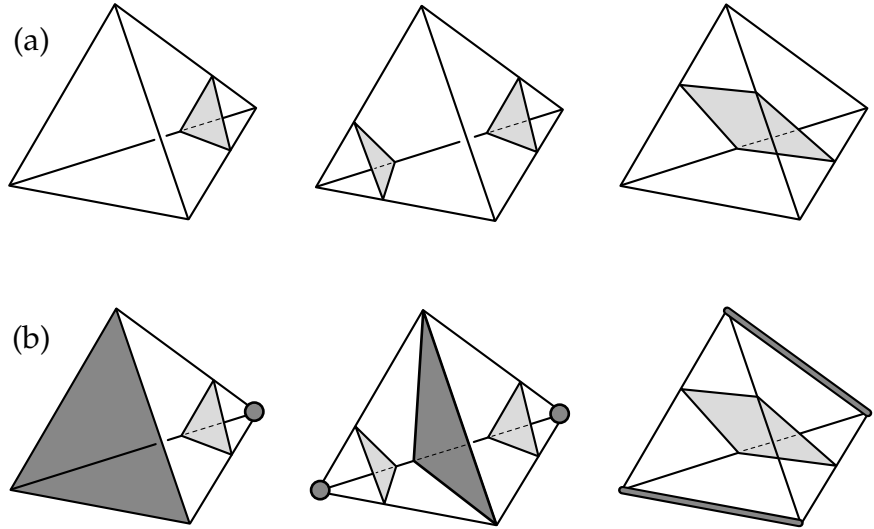

\centering\leavevmode
\begin{xyoverpic*}{(408,249)}{scale=0.8}{pictures/normal}
,(11,230)*{\mbox{(a)}}
,(11,89)*{\mbox{(b)}}
\end{xyoverpic*}
  \caption{
    At the top, the three possibilities, besides the empty set, for how
    $\Sigma$ intersects a tetrahedron of $\mathcal T$.  The pictures at the
    bottom describe, in dark grey, a 2-complex $P$ to which $N \setminus \Sigma$ deformation
    retracts.}
  \label{fig_simple_spine}
\end{figure}

To check this, we consider the 8-fold cover $N$ of the orbifold $B$ of
Figure~\ref{B-pic} given by
\[
\pi_1(B) \to S_8 \quad \mbox{where} \quad a \mapsto (1, 2, 3, 6, 8,
5, 7, 4), b \mapsto (1, 3)(2, 5)(4, 6)(7, 8),
\]
which is a manifold.  Now $\pi_1(N)$ is not a subgroup of $\pi_1(M)$,
but we can consider their minimal common cover $Y$ with fundamental
group $\pi_1(N) \cap \pi_1(M)$ as shown
\[
\xymatrix @ur {
N \ar[d] & Y  \ar[l] \ar[d] \\
B & M  \ar[l]
}
\]
A key fact is that $\dim H^1(Y; \R) = 3$.  Thus we may apply the
theorem of Stallings mentioned in Section~\ref{sec_norm_cover} and
conclude, as $Y$ and $M$ have the same first Betti number, that $M$
fibers if and only if $Y$ does. Thus if $N$ fibers, so does $M$. Now
$H_1(N ; \Z) = \Z/4 \oplus \Z/8 \oplus \Z$, and its Alexander
polynomial is $\Delta_N = t^4 + 30 t^2 + 1$, suggesting that $N$ is
fibered by genus 2 surfaces. The method of
Section~\ref{subsec_ex_thurston} finds a genus 2 normal surface
$\Sigma$ representing the generator of $H_2(N; \Z)$ (the particular
triangulation of $N$ we used, and hence an explicit description of
$\Sigma$, is available at \cite{paperwebsite}). To show that $N
\setminus \Sigma$ is $\Sigma \times I$ and hence $N$ fibers, it
suffices to check that $\pi_1(N \setminus \Sigma)$ is
\emph{abstractly} isomorphic to $\pi_1(\Sigma)$ (see
e.g.~\cite[Thm.~10.6]{HempelBook}).  To find an initial presentation
of $\pi_1(N \setminus \Sigma)$, we can exploit the fact that in a
suitable triangulation $\mathcal T$ of $N$, the normal surface
$\Sigma$ is very simple; in each tetrahedron it looks like one of
the three possibilities shown in Figure~\ref{fig_simple_spine}(a).
Part (b) of the same figure describes a 2-complex $P$ which is a spine
for $N \setminus \Sigma$, i.e.~$N \setminus \Sigma$ deformation retracts to $P$.  It is
easy to read off a presentation for $\pi_1(N \setminus \Sigma)$ from this spine;
simplifying this presentation using Tietze transformations
\cite{MagnusKarrassSolitar} and relabeling the generators yields
\[
\pi_1(N \setminus \Sigma) = \spandef{a,b,c,d}{a b a^{-1} b^{-1} c d c^{-1} d^{-1} = 1}
\]
which is the standard presentation for the fundamental group of a genus
2 surface.  So $N$ fibers, and hence so does $M$.

\begin{remark}
  This procedure for checking fibering is applicable much more
  generally.  For any normal surface, there is
  a natural spine for its complement, introduced by Casson in his
  study of 0-efficient triangulations (paralleling
  \cite{JacoRubinstein0eff}).  Writing down a presentation of $\pi_1(N
  \setminus \Sigma)$ from this is straightforward.  In more complicated cases,
  simplifying the presentation of $\pi_1(N \setminus \Sigma)$ in the hopes of
  recognizing it as a surface group, it is important to employ not
  just Tietze transformations, but Nielsen transformations as well,
  i.e.~the application of elements in $\Aut(F_n)$ to the generating
  set of the presentation in order to shorten the relators.  If one
  can reduce down to a presentation with only one relator, then
  checking it is a surface group is often easy: just glue up the edges
  of the relator polygon, and if there is only one vertex in what
  results, then the group is indeed a surface group.
\end{remark}

\subsection{Oldforms for $M = X(\q\qbar)$}

We turn now to the last assertion of Theorem~\ref{thm_top_prop},
namely that the space of oldforms in $H^1(M; \Q)$ is the line
determined by the barycenters of the non-fibered faces.  By
Theorem~\ref{thm_autom_prop}, we know that $\dim H^1(M;\Q)\newforms \geq
2$, so as $\dim H^1(M; \Q) = 3$ it follows that $H^1(M; \Q)\oldforms$
has dimension at most $1$.  Since $H^1(X(\q); \Q)$ is nonzero, it
follows that $H^1(M; \Q)\oldforms$ is
the image of $\phi^* \maps H^1(X(\q); \Q) \to H^1(M;\Q)$, where
$\phi = \phi_1 \maps M \to X(\q)$ is the standard covering map.
Since we have explicit presentations for both $\pi_1(X(\q))$ and
$\pi_1(M)$ as subgroups of $\pi_1(B)$, computing the image of
$\phi^*$ is straightforward; for details see \cite{paperwebsite}.

\section{Proof of the main result}\label{sec_proof_of_main_thm}

\subsection{Proof of Theorem~\ref{rapid_growth}}  Continuing the notation of
Section~\ref{sec_example}, let $M$ be the arithmetic hyperbolic
\3-manifold described there, and $\omega \in H^1(M; \Z)$ be the
class given by part (\ref{ex_top_part}) of
Theorem~\ref{thm_example_properties}. Then $\omega$ is fibered and
for each prime $\p$ in the special set $\sP_K$ we have $T_\p(\omega)
= 0$.  Denote the primes of $\sP$ as $p_1, p_2, p_3, \ldots$ in
increasing order, and fix an ordering, once and for all, of the
primes of $\sP_K$ as $\p_1,\p_2, \ldots$ such that $i\leq j \,
\implies \, N(P_i)\leq N(P_j)$.

Consider the congruence cover $M_n = M(\p_1 \cdots \p_{2n})$.  By
Theorem~\ref{cong_cover_thm}, the number $\nu_n$ of pairs of fibered
faces of $M_n$ is at least $2^{2n}$, and the degree $d_n$ of $M_n \to
M$ is (see Section~\ref{cong_covers}):
\begin{equation}\label{prod_form_of_dn}
d_n = \prod_{i=1}^{2n} \left(1 + \Norm_{K/\Q}(\p_i)\right) = \left(\prod_{i=1}^n\left(1 + p_i\right)\right)^2
\end{equation}
To prove Theorem~\ref{rapid_growth} we will estimate $d_n$ by using
the fact that $\sP$ consists of $1/4$ of all rational primes.  More
precisely, recall that the \emph{natural density} of a set
$\mathcal{Y}$ of primes is $\alpha\in [0,1]$ if the number of primes
in $\mathcal{Y}$ which are $\leq x$ is asymptotic (for $x$ large) to
$\alpha$ times $\pi(x)=\vert \{p\leq x\}\vert$.  In our case, $\sP$
is the set of primes $p$ which split in both $\Q[\sqrt{-3}]$ and
$\Q[\sqrt{-7}]$, i.e., those which split completely in the
biquadratic field $\Q[\sqrt{-3},\sqrt{-7}]$, and so $\sP$ has
natural density $1/4$ by the Tchebotarev density theorem.
Theorem~\ref{rapid_growth} now follows directly from:
\begin{lemma}\label{lemma_bound}
  Let $\sP = \{p_1, p_2, p_3 \ldots \}$ be a set of rational primes of
  natural density $1/4$, and $d_n$ the product defined in
  (\ref{prod_form_of_dn}).  Then
  \[
  n =  \left(\frac12 +  o(1)\right)\frac{\log d_n}{\log\log d_n}.
  \]
  In particular,
  \[
  2^{2n} = \exp\left(\left(\log 2+ o(1)\right)\frac{\log d_n}{\log\log
      d_n}\right) \, \leq c_t e^{(\log d_n)^t},
  \]
  for any $t<1$, for a suitable constant $c_t>0$.
\end{lemma}
Some people prefer to work with weaker notions of density, like
analytic density, but in those cases this lemma may not hold.

\begin{proof}
  If $f$ and $g$ are functions of $x$ on a subset of $\R_+$, we will
  write $f(x) \sim g(x)$ if their ratio tends to $1$ as $x$ goes to
  infinity.  Now recall that the prime number theorem asserts that
  $\pi(x) \sim \frac{x}{\log x}$, and consequently the $n^{\mathrm{th}}$ prime is
  $\sim n\log n$. As the set $\sP$ has natural density $1/4$,
  the $n^\mathrm{th}$ prime $p_n$ in $\sP$ satisfies
  \begin{equation}\label{eqn_pn}
  p_n  \sim  4n\log n.
  \end{equation}
  Moreover, we claim that
  \begin{equation}\label{eqn_log_dn}
  \log\left(\prod\limits_{j=1}^n (1+p_j)\right) \sim \frac14  p_n.
  \end{equation}
  In fact, the Riemann hypothesis implies that this holds with an
  additive error term $\rho(n)$ of the order of $n^{1/2+\varepsilon}$.  Here, we
  only need that $\rho(n)$ is $o(n\log n)$, which is known; in
  fact one has (cf. \cite[Ch.~6]{MontgomeryVaughan}):
  \[
  \abs{\rho(n)} \leq ne^{-c\sqrt{\log n}} \quad \mbox{for some $c>0$.}
  \]

  Combining the asymptotic statements (\ref{eqn_pn}) and
  (\ref{eqn_log_dn}), and noting the formula for $d_n$, we get
  \[
  2n\log n \sim \log d_n.
  \]
  This yields, upon taking logarithms,
  \[
  \log n +\log(2\log n)  \sim  \log\log d_n,
  \]
  and since the left hand side is $\left(1+o(1)\right)\log n$, we obtain
  \[
  2n \sim \frac{\log d_n}{\log n} = \left(1+o(1)\right)\frac{\log d_n}{\log\log d_n}.
  \]
  Since $2^{2n}= e^{2n\log 2}$, this implies
  \begin{equation}\label{eqn_2n}
  2^{2n} = {\rm exp}\left(\left(\log 2+o(1)\right)\frac{\log d_n}{\log\log
  d_n}\right).
  \end{equation}
  Moreover, for any real number $t<1$, the quantity $\frac{\log
    d_n}{\log\log d_n}$ dominates $(\log d_n)^t$ as $n \to \infty$. So the
  right hand side of (\ref{eqn_2n}) is bounded below by a constant
  (depending on $t$) times ${\rm exp}\left((\log d_n)^t\right)$, as
  asserted by the lemma.  This completes the proof of
  Theorem~\ref{rapid_growth}.
\end{proof}

\subsection{A lower bound for $b_1(M_n)$}
Recall that the first Betti number of any of the manifolds $M_n$ is
the same as the dimension of space of the cohomological cusp forms
on $D^\ast$ of level $I_n = \q \qbar \prod\limits_{j=1}^n
\p_j{\overline \p}_j$, with $\p_j, {\overline \p}_j$ being the prime
ideals in $\oO_K$ above the rational prime $p_j\in\sP$. From the
arithmetic point of view, the dimension of this space is a complete
mystery. This is because the two standard tools, namely the
Riemann-Roch theorem and the Selberg trace formula, which work so
well in the case of Hilbert modular forms (of weight $2$), are not
applicable here; the former because the $M_n$ are not algebraic
varieties, and the latter because the relevant archimedean
representation is not in the discrete series, making it impossible
to separate it from all the other archimedean representations.
Nevertheless, we can get a lower bound as follows.

First recall that for any $N\geq 1$, the dimension of the space $S_2(N)$
of classical, holomorphic cusp forms on the upper half-plane $\mathcal
H$ of level $N$ is the genus $g_0(N)$ of the standard (cusp)
compactification of the Riemann surface $\Gamma_0(N)\backslash {\mathcal H}$,
where $\Gamma_0(N)$ is the subgroup of $\SL{2}{\Z}$ consisting of matrices
$\mysmallmatrix abcd$ with $N \mid c$. It is classical
that for $N$ square-free, one has
\[
g_0(N) = 1+\frac{1}{12} \prod\limits_{p \mid N} (p+1)
-\frac{2^s}{4}-\frac{2^t}{3}-\frac{2^m}{2},
\]
where $m$ is the number of prime divisors of $N$, and $s$
(resp.~$t$) the number of such divisors which are $1$ mod $4$
(resp.~$1$ mod $3$); the last three terms of the formula are the error
terms corresponding to the elliptic points $i$ and $(1+\sqrt{-3})/2$
on $\mathcal H$ and the cusp at infinity. Applying this with $N=N_n=
7\prod\limits_{j=1}^n p_j$, with $p_j \neq 7$ being split in
$\Q[\sqrt{-3}]$ and inert in $\Q[\sqrt{-7}]$, we see by
the formula for $d_n$ that
\[
g_0(N_n) \geq \frac{\sqrt{d_n}-(13)2^n}{12}.
\]
On the other hand, by Lemma~\ref{lemma_bound}, the term $2^n$ is
bounded from above by any positive power of $d_n$. Moreover, the
dimension of the subspace $S_2(N_n)^{7\mbox{-}\mathrm{new}}$ of $S_2(N_n)$
consisting of cusp forms of level $N$ which are new at $7$ is $\frac57
g_0(N)$. Putting these together, we obtain, for any $\varepsilon > 0$,
\[
\dim_\C S_2(N_n)^{7\mbox{-}\mathrm{new}} \geq
\frac{5}{84}d_n^{\frac12-\varepsilon},
\]
for all large enough $n$ (depending on $\varepsilon$). It is an easy
exercise to see that the dimension on the left is derived purely
from the knowledge of the dimensions of the space
$S_2(N_n')^{\mathrm{new}}$ of {\it newforms} of level $N_n'$ with $7
\mid N_n' \mid N_n$, together with the factorization of $N_n/N_n'$.

Next consider the base change $g \to g_K$ from $\GLtwo/\Q$ to
$\GLtwo/K$ (\cite{Langlands}). We claim that the image of each
$S_2(N_n')^{\mathrm{new}}$ under this base change has the same
dimension. Indeed, if $g_K=h_K$ for two newforms $g, h$ in this
space, the newform $h$ must be a twist of $g$ by the quadratic Dirichlet
character $\delta$ corresponding to $K$. But since $g, h$ have odd
levels, this would imply that $h$ has level divisible by $4$, the
level of $\delta$, which is impossible. Hence the claim. Now by
construction, every newform $g_k$, or rather the cuspidal
automorphic representation $\beta_K$ of $\GL{2}{\A_K}$ attached to
it, is Steinberg at each of the primes dividing $7$, and so
transfers to our $D^\ast$ and defines a cusp form $\beta_K^D$
relative to $\Gamma_D(N_n/7)$; it also lifts to level
$\Gamma_D(N_n)$ by the natural map.  The cusp forms on $\GLtwo/K$
which are Steinberg at $\{ \q, \qbar \}$, and which are old away from $7$,
are seen to correspond in a one-to-one way to old forms on
$D^\ast$, since $D_v^\ast\simeq \GL{2}{K_v}$ for any $v\nmid 7$.
Consequently, we get the following
\begin{proposition} \label{prop_lower_bound}
  Let $\{M_n\}$ be the arithmetic tower of closed hyperbolic
  $3$-manifolds we have constructed above, with $M_n$ a non-normal
  cover of $M$ of degree $d_n$. Then for every $\varepsilon
  >0$, we can find an integer $L_\varepsilon$ such that
\[
n\geq L_\varepsilon \, \implies \, b_1(M_n) \geq
\frac{5}{84}d_n^{\frac12-\varepsilon}.
\]
\end{proposition}

\section{The Whitehead link complement} \label{sec-whitehead}

\begin{figure}[hbt]
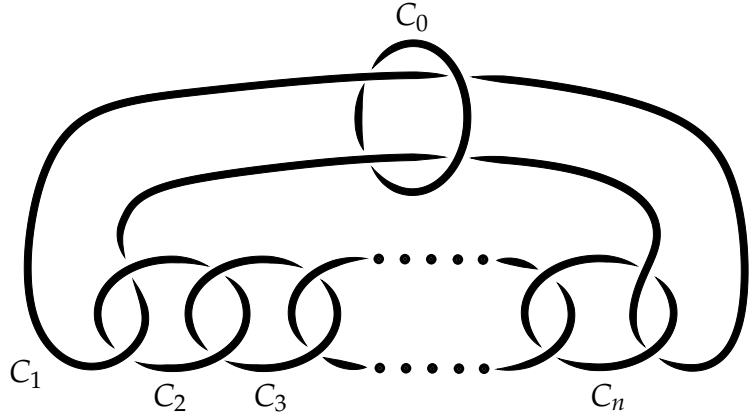

\centering\leavevmode
\begin{xyoverpic*}{(303,142)}{scale=0.9}{pictures/Ln}
,(162,141)*+!D{C_0}
,(12,10)*+!UR{C_1}
,(61,4)*++!U{C_2}
,(103,4)*++!U{C_3}
,(244,4)*++!U{C_n}
\end{xyoverpic*}
\caption{The covering manifold $W_n$ is the exterior of the shown link $L_n$.}
\label{fig-n-cover}
\end{figure}

In this section, we give a very concrete example of covers of a finite
volume hyperbolic 3-manifold with cusps where the number of fibered
faces increases exponentially in the degree.
\begin{figure}[h]
\includegraphics[scale=0.6]{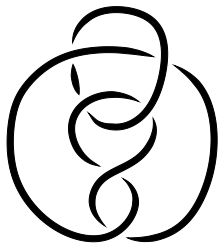}
\caption{The Whitehead link.}
\label{fig-whitehead}
\end{figure}
Consider the Whitehead link shown in Figure~\ref{fig-whitehead}, and
let $W$ be its exterior.  The interior of $W$ admits a hyperbolic
metric of finite volume; indeed, $W$ is arithmetic with $\pi_1(W)$ a
subgroup of $\PSL{2}{\Z[i]}$ of index 12.  We consider $n$-fold cyclic
covers dual to the thrice punctured sphere bounded by one of the
components.  More precisely, let $\omega$ be an element of the basis of
$H^1(W; \Z)$ which is dual to a basis of $H_1(W ; \Z) \cong \Z^2$
consisting of meridians.  (This link is quite symmetric, so
it doesn't matter which basis element we pick.)  Then let $W_n$ be
the cover corresponding to the kernel of the composition $\pi_1(W)
\stackrel{\omega}{\to} \Z \to \Z/n$.  The manifold $W_n$ is also the
complement of the link $L_n$ pictured in Figure~\ref{fig-n-cover}.
Following that figure, let $\mu_i \in H_1(W_n ; \Z)$ be a
meridian for the component $C_i$ of $L_n$.  Let $\{ \omega_i \}$ be the
dual basis for $H^1(W_n; \Z)$.  This section is devoted to proving:
\begin{theorem}\label{thm_whitehead}
  Let $W_n$ be the $n$-fold cyclic cover of the exterior of the
  Whitehead link described above.  Then the number of fibered faces of
  the Thurston norm ball $B_T$ of $W_n$ is $2^{n+1}$.  More precisely,
  $B_T$ is the convex hull of the $\omega_i$ and all of its top-dimensional faces
  are fibered.
\end{theorem}
The polytope $B_T$ is called an $n$-orthoplex or a cross polytope;
it is dual to the unit cube in $H_1(W_n; \Z)$ spanned by the
$\mu_i$. Note that this form for $B_T$ is plausible as each
component of $L_n$ bounds an obvious twice-punctured disc, and hence
$\norm{\omega_i}_T = 1$ for every $i$ as hyperbolic manifolds
contain no simpler essential surfaces.   The theorem will follow
easily from
\begin{lemma}\label{lem_whitehead}
  Any class of $H^1(W_n; \Z)$ of the form
  \[
  \omega = \epsilon_0 \omega_0 + \epsilon_1 \omega_1 + \cdots + \epsilon_n \omega_n  \quad \mbox{for $\epsilon_i \in \{ -1, 1 \}$}
  \]
  fibers and $\norm{\omega}_T = n + 1$.
\end{lemma}

Let us first derive the theorem from the lemma.

\begin{proof}[Proof of Theorem~\ref{thm_whitehead}]
  We need to show that for any choices of $\epsilon_i \in \{ -1, 1\}$, the
  simplex $\Delta$ spanned by $\eta_i = \epsilon_i \omega_i$ is a face of $B_T$.  We
  know $\norm{\eta_i}_T = 1$, and so $\Delta \subset B_T$.  By
  Lemma~\ref{lem_whitehead}, the class $\eta = \frac{1}{n + 1} \sum \eta_i$
  has norm $1$; this gives another point in $\Delta$ which we know lies in
  $\partial B_T$.  Now consider $\alpha \in \Delta$, which necessarily has the form
  \[
  \alpha = \sum_{i = 0}^n a_i \eta_i \quad \mbox{for $a_i \geq 0$ and $\sum a_i = 1$.}
  \]
  We need to show $1 \leq \norm{\alpha}_T$, as then $\alpha$ lies in $\partial B_T$;
  knowing this for all $\alpha$ implies $\Delta$ is a face of $B_T$.

 After permuting the coordinates, we can
  assume $a_0 \geq a_i$ for all $i$.  As a consequence, $a_0 \geq 1/(n +
  1)$ and
  \[
  \eta = b_0 \alpha + b_1 \eta_1 + b_2 \eta_2 + \cdots + b_n \eta_n \quad \mbox {where $b_i \geq 0$ and $\sum_{i=0}^n b_i = 1$.}
  \]
  Now the sublinearity of the Thurston norm gives
  \[
  \norm{\eta}_T \leq b_0 \norm{\alpha}_T + b_1 \norm{\eta_1}_T + b_2 \norm{\eta_2}_T + \cdots + b_n \norm{\eta_n}_T.
  \]
  Evaluating the norms we know gives:
  \[
  1 \leq b_0 \norm{\alpha}_T + b_1 + \cdots + b_n
  \]
  Since $\sum b_i = 1$, this implies $1 \leq \norm{\alpha}_T$, as required
  to prove the theorem.
\end{proof}

We conclude this section with

\begin{proof}[Proof of Lemma~\ref{lem_whitehead}]
  Since the Thurston norm is invariant under $\omega \mapsto -\omega$, we can
  assume $\epsilon_0 = 1$.  We analyze this class using Murasugi sum,
  closely following the approach used in \cite{Leininger2002} to
  analyze chain links.
  \begin{figure}
    \includegraphics[scale=0.9]{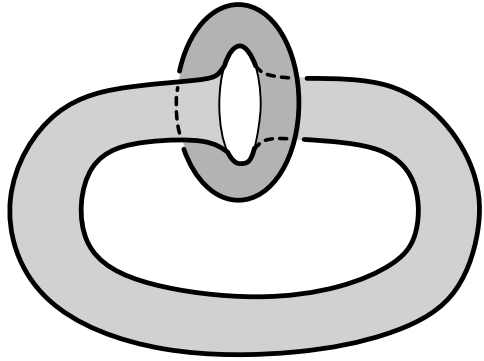}
    \caption{The base link/surface for the Murasugi sum shown in Figure~\ref{fig-link-lai}.}
    \label{fig-murasugi-base}
  \end{figure}
  Given $a_i \in \{ -1, 1\}$, consider the link $L(a_1, a_2, \ldots, a_n)$
  shown in Figure~\ref{fig-link-lai}, which is described as the
  boundary of a surface which is the Murasugi sum of the surface in
  Figure~\ref{fig-murasugi-base} with $n$ Hopf bands, where the twist
  on the $i^{\mathrm{th}}$ Hopf band is right- or left-handed
  depending on $a_i$.  Let $M(a_1,\ldots, a_n)$ be the exterior of $L(a_1,
  \ldots, a_n)$.  Since all the surfaces in the Murasugi sum are fibers,
  the surface $\Sigma$ pictured in Figure~\ref{fig-murasugi-base} is a
  fiber of a fibration of $M(a_1,\ldots, a_n)$ over the circle
  \cite{Stallings1978, Gabai1985sum}.
   \begin{figure}

     \includegraphics[scale=0.9]{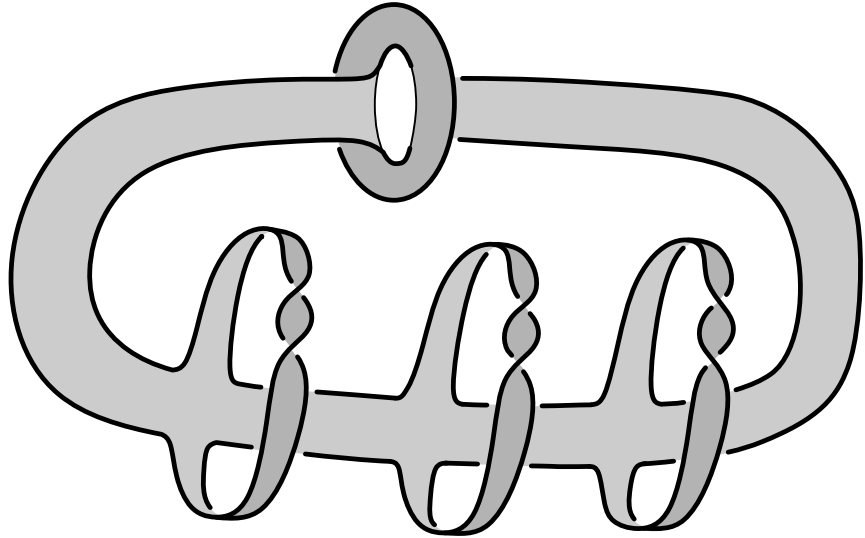}
    \caption{
      The link $L(a_1, a_2, \ldots, a_n)$ is the Murasugi sum of the link
      in Figure~\ref{fig-murasugi-base} with $n$ Hopf bands with
      appropriately oriented twists.    Shown is $L(1, 1, -1)$.}
    \label{fig-link-lai}
  \end{figure}
  Orient the components of $L(a_1, \ldots, a_n)$ as the boundary of the
  surface shown in Figure~\ref{fig-link-lai}, and let $\omega_i$ be the
  corresponding dual basis of $H^1(M(a_1,\ldots, a_n) ; \Z)$.  Thus $\Sigma$
  is Poincar\'e dual to $\omega = \omega_1 + \cdots + \omega_n$ and $\norm{\omega}_T = n +
  1$.  The key to the lemma is to show:
  \begin{claim}\label{whitehead-claim}
    Fix $a_1, \ldots, a_n$ in $\{-1, 1\}$, and pick $k \in \{1, 2, \ldots,
    n\}$.  Suppose $b_i = a_i$ except $b_{k - 1} = -a_{k - 1}$ and
    $b_k = -a_k$. (If $k = 1$ then here $k - 1$ should be interpreted
    as $n$.)  Then $M(a_1,\ldots,a_n) \cong M(b_1, \ldots, b_n)$ via a
    homeomorphism that acts on $H^1$ by
    \begin{align*}
      \omega_i \mapsto \omega_i &\quad \mbox{for $i \neq k$ and } \\
      \omega_k \mapsto -\omega_k&
    \end{align*}
  \end{claim}
  The proof of this claim is given in Figure~\ref{fig-proof-of-claim}.
  \begin{figure}
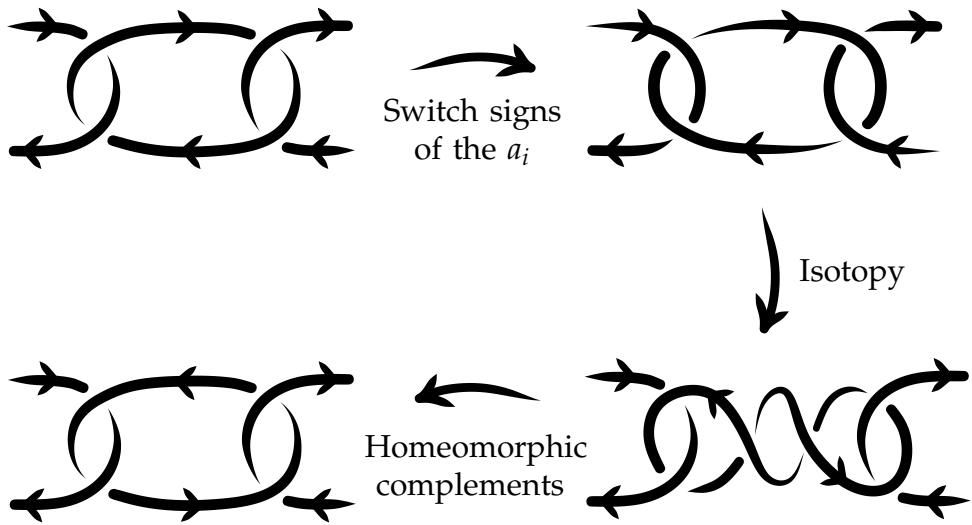

    \centering\leavevmode
    \begin{xyoverpic*}{(365,196)}{scale=1.0}{pictures/pf-of-claim}
      ,(176,175)*++++!U{\txt<2.5cm>{Switch signs of  the $a_i$}}
      ,(289,94)*+++!L{\txt{Isotopy}}
      ,(175,51)*+++++!U{\txt<2.8cm>{Homeomorphic complements}}
    \end{xyoverpic*}
    \caption{
      Switching the signs of two consecutive $a_i$ can be partially
      reversed at the cost of adding a twist as shown in the bottom
      right picture.  (Here both $a_i$ are $1$, but the other cases
      are similar.)  We can undo said twist by a homeomorphism of the
      link complements by cutting along a twice punctured disc bounded
      by component $C_0$ of the link and doing a full twist. To prove
      Claim~\ref{whitehead-claim}, just note the different
      orientations of the central loop in the two pictures at left.  }
    \label{fig-proof-of-claim}
  \end{figure}
  Now to prove the lemma, note that $L(1,1, \ldots, 1)$ is $L_n$ and
  applying the claim repeatedly to a suitable starting $L(a_1,\ldots,a_n)$
  allows us to see any that class of the form $\omega_0 + \epsilon_1 \omega_1 + \cdots + \epsilon_n
  \omega_n$ is represented by a fiber surface with Euler characteristic
  $-(n + 1)$ as required.
\end{proof}

\section{Work of Long and Reid}
\label{sec_long_reid}

In this section, we give a proof of Theorem~\ref{thm_long_reid} which
is simpler than the original one in \cite{LongReid2007}.  A different
simplification was given by Agol in \cite{Agol2007}. We begin by
introducing an invariant of a fibered cohomology class $\omega$ of a
hyperbolic manifold $M$.  Let $\F_\omega$ be the associated transverse
pseudo-Anosov flow described in Section~\ref{subsec_frieds_work}.  As
discussed, this flow has a positive finite number of singular orbits,
each of which is closed.  From the closed geodesics $\gamma_1, \gamma_2, \ldots,
\gamma_n$ homotopic to these orbits, define $c(\omega)$ to be the subset of
the universal cover $\hthree$ of $M$ consisting of the inverse images
of the $\gamma_i$.  By Theorem~\ref{thm_fried}, any other class $\eta$ in
the fibered face of $\omega$ has isotopic $\F_\eta$, and hence $c$ is really
an invariant of a fibered face.  In fact, we regard the geodesics in
$c(\omega)$ as unoriented, so that $c$ is an invariant of a fibered face
pair.  With this invariant in hand, we turn to the proof itself.

\begin{proof}[Proof of Theorem~\ref{thm_long_reid}]
  Let $\omega$ be a fibered class for our arithmetic 3-manifold $M = \Gamma \backslash
  \hthree$.  We will first find a cover of $M$ with two pairs of
  fibered faces.  Following Section~\ref{sec-source-hecke}, given $g
  \in \Comm(\Gamma)$ we have the associated manifold $M_g$ with two
  covering maps $p_g, q_g \maps M_g \to M$.  Then $M_g$ has two fibered
  classes $p_g^*(\omega)$ and $q_g^*(\omega)$, with associated $c$ invariants
  equal to $c(\omega)$ and $g^{-1} \cdot c(\omega)$ respectively.  Thus these two
  fibered classes of $M_g$ lie in genuinely distinct faces provided
  $c(\omega)$ is not setwise invariant under $g$. Now the setwise
  stabilizer of $c(\omega)$ in $\PSL{2}{\C}$ is a discrete subgroup which
  contains $\Gamma$ as a subgroup of finite index; as $\Comm(\Gamma)$ is dense
  in $\PSL{2}{\C}$, there are plenty of $g \in \Comm(\Gamma)$ which do not
  fix $c(\omega)$, as desired.

  To create any number of fibered faces, we simply repeat this process
  inductively, assuming at each stage we have constructed a cover
  $M_n$ with fibered classes $\omega_1,\ldots, \omega_n$ where the $c(\omega_i)$ are
  all distinct and then building a cover of using a $g \in \Comm\left(\pi_1(M_n)\right)$
  so that all $c(\omega_i)$ and $g \cdot c(\omega_k)$ are distinct.
\end{proof}

{\RaggedRight \bibliographystyle{math} \bibliography{hecke_faces} }

\end{document}